\documentclass{article}
\usepackage{amsmath}
\usepackage{amssymb}
\usepackage{geometry}
\usepackage{titlesec}
\usepackage[hyperfootnotes=false,linktoc=none]{hyperref}
\hypersetup{
    colorlinks=true,
    linkcolor=black,
    urlcolor=blue,
    citecolor=blue
}

\titleformat{\section}
  {\normalfont\Large\bfseries\centering}{\thesection}{1em}{}
\titlespacing*{\section}{0pt}{3.5ex plus 1ex minus .2ex}{2.3ex plus .2ex}

\newcommand{\tnref}[1]{\hyperlink{#1}{[#1]}}

\title{{\small\textit{English translation of Sophie Kowalevski's}} \\[0.5em] \MakeUppercase{On the Problem of the Rotation \\ of a Rigid Body About a Fixed Point}}
\author{Sophie Kowalevski\footnote{This memoir is the summary of a work to which the Academy of Sciences of Paris, in its solemn session of December 24, 1888, awarded the Bordin prize, raised from 3000 to 5000 francs.} \\ Stockholm}
\date{}

\begin{document}

\maketitle

\section*{§ 1.}

The problem of the rotation of a heavy rigid body about a fixed point can be reduced, as is well known, to the integration of the following system of differential equations:

\begin{equation}
\begin{aligned}
A\frac{dp}{dt} &= (B - C)qr + Mg(y_0\gamma'' - z_0\gamma'), & \frac{d\gamma}{dt} &= r\gamma' - q\gamma'', \\
B\frac{dq}{dt} &= (C - A)rp + Mg(z_0\gamma - x_0\gamma''), & \frac{d\gamma'}{dt} &= p\gamma'' - r\gamma, \\
C\frac{dr}{dt} &= (A - B)pq + Mg(x_0\gamma' - y_0\gamma), & \frac{d\gamma''}{dt} &= q\gamma - p\gamma'.
\end{aligned}
\end{equation}
The constants $A$, $B$, $C$, $Mg$, $x_0$, $y_0$, $z_0$ that appear in these equations have the following meaning:
\begin{itemize}
\item $A$, $B$, $C$ are the principal axes of the ellipsoid of inertia of the considered body, with respect to the fixed point;
\item $M$ is the mass of the body;
\item $g$ is the intensity of the gravitational force;
\item $x_0$, $y_0$, $z_0$ are the coordinates of the center of gravity of the considered body in a coordinate system whose origin is at the fixed point and whose direction coincides with that of the principal axes of the ellipsoid of inertia.
\end{itemize}
Until now, one has only succeeded in integrating these equations in two particular cases:
\begin{enumerate}
\item The case of Poisson (or of Euler) where one has $x_0 = y_0 = z_0 = 0$,
\item The case of Lagrange where one has $A = B$, $x_0 = y_0 = 0$.
\end{enumerate}
In these two cases, the integration is performed with the aid of elliptic functions whose argument is a linear entire function of time.
The six quantities $p$, $q$, $r$, $\gamma$, $\gamma'$, $\gamma''$ are in these two cases uniform functions of time, having no other singularities than poles for all finite values of the variable.

Do the integrals of the considered differential equations preserve this property in the general case?
If so, one would need to be able to integrate these differential equations with the aid of series of the form

\begin{equation}
\begin{aligned}
p &= t^{-n_1}(p_0 + p_1t + p_2t^2 + \cdots), \\
q &= t^{-n_2}(q_0 + q_1t + q_2t^2 + \cdots), \\
r &= t^{-n_3}(r_0 + r_1t + r_2t^2 + \cdots), \\
\gamma &= t^{-m_1}(f_0 + f_1t + f_2t^2 + \cdots), \\
\gamma' &= t^{-m_2}(g_0 + g_1t + g_2t^2 + \cdots), \\
\gamma'' &= t^{-m_3}(h_0 + h_1t + h_2t^2 + \cdots),
\end{aligned}
\end{equation}
where $n_1$, $n_2$, $n_3$, $m_1$, $m_2$, $m_3$ denote positive integers, and these series, in order to represent the general system of integrals of the considered differential equations, should contain five arbitrary constants.
It is therefore necessary to examine whether such an integration is possible.
One easily verifies, by comparing the exponents of the first terms in the left and right-hand sides of the considered equations, that one must have
\[
n_1 = n_2 = n_3 = 1, \quad m_1 = m_2 = m_3 = 2.
\]
One then finds, by setting for brevity:
\[
A_1 = B - C, \quad B_1 = C - A, \quad C_1 = A - B
\]
and by choosing the unit of length such that $Mg = 1$, that the coefficients $p_0$, $q_0$, $r_0$ must satisfy the following equations:

\begin{equation}
\begin{aligned}
-Ap_0 &= A_1q_0r_0 + y_0h_0 - z_0g_0, & -2f_0 &= r_0g_0 - q_0h_0, \\
-Bq_0 &= B_1r_0p_0 + z_0f_0 - x_0h_0, & -2g_0 &= p_0h_0 - r_0f_0, \\
-Cr_0 &= C_1p_0q_0 + x_0g_0 - y_0f_0, & -2h_0 &= q_0f_0 - p_0g_0.
\end{aligned}
\end{equation}
For the last three of these equations to be satisfied by non-zero values of $f_0$, $g_0$, $h_0$, it is necessary that the following determinant is zero:
\[
\begin{vmatrix}
2 & -r_0 & q_0 \\
r_0 & 2 & -p_0 \\
-q_0 & p_0 & 2
\end{vmatrix} = 2(4 + p_0^2 + q_0^2 + r_0^2).
\]
If one adds the first three of equations (3) after having multiplied them respectively by $Ap_0$, $Bq_0$, $Cr_0$, one finds
\begin{equation}
A^2p_0^2 + B^2q_0^2 + C^2r_0^2 = x_0(Bq_0h_0 - Cr_0g_0) + y_0(Cr_0f_0 - Ap_0h_0) + z_0(Ap_0g_0 - Aq_0f_0);
\end{equation}
but from the last three of equations (3) it follows that

\begin{equation}
\begin{aligned}
2(Bq_0h_0 - Cr_0g_0) &= p_0(Ap_0f_0 + Bq_0g_0 + Cr_0h_0) - f_0(Ap_0^2 + Bq_0^2 + Cr_0^2), \\
2(Cr_0f_0 - Ap_0h_0) &= q_0(Ap_0f_0 + Bq_0g_0 + Cr_0h_0) - g_0(Ap_0^2 + Bq_0^2 + Cr_0^2), \\
2(Ap_0g_0 - Bq_0f_0) &= r_0(Ap_0f_0 + Bq_0g_0 + Cr_0h_0) - h_0(Ap_0^2 + Bq_0^2 + Cr_0^2).
\end{aligned}
\end{equation}
The six quantities $p$, $q$, $r$, $\gamma$, $\gamma'$, $\gamma''$ must satisfy the two algebraic relations
\begin{align*}
Ap^2 + Bq^2 + Cr^2 &= 2(x_0\gamma + y_0\gamma' + z_0\gamma'') + l_1, \\
Ax_0p + By_0q + Cz_0r &= l,
\end{align*}
where $l$ and $l_1$ denote integration constants, so it is necessary that
\begin{align*}
Ap_0f_0 + Bq_0g_0 + Cr_0h_0 &= 0, \\
x_0f_0 + y_0g_0 + z_0h_0 &= \frac{1}{2}(Ap_0^2 + Bq_0^2 + Cr_0^2).
\end{align*}
It follows therefore from equations (4) and (5) that
\[
A^2p_0^2 + B^2q_0^2 + C^2r_0^2 = -\frac{1}{4}(Ap_0^2 + Bq_0^2 + Cr_0^2)^2.
\]
It is now necessary to distinguish two cases.

\textbf{First case.} Suppose that none of the following quantities are zero
\[
A_1 = B - C, \quad B_1 = C - A, \quad C_1 = A - B.
\]
If one then sets
\begin{align*}
\frac{1}{2}(p_0^2 + q_0^2 + r_0^2) &= \lambda_0, \\
\frac{1}{2}(Ap_0^2 + Bq_0^2 + Cr_0^2) &= \lambda, \\
\frac{1}{2}(A^2p_0^2 + B^2q_0^2 + C^2r_0^2) &= \lambda_1,
\end{align*}
one has
\begin{align*}
\frac{1}{2}p_0^2 &= -\frac{BC\lambda_0 - (B + C)\lambda + \lambda_1}{B_1C_1}, \\
\frac{1}{2}q_0^2 &= -\frac{CA\lambda_0 - (C + A)\lambda + \lambda_1}{C_1A_1}, \\
\frac{1}{2}r_0^2 &= -\frac{AB\lambda_0 - (A + B)\lambda + \lambda_1}{A_1B_1}.
\end{align*}
Now we have found
\[
\lambda_0 = -2, \qquad \lambda_1 = -\frac{1}{2}\lambda^2.
\]
One therefore has
\begin{equation}
\begin{aligned}
p_0^2 &= \frac{4BC\lambda_0 - (B + C)^2\lambda + \lambda^2}{B_1C_1} = \frac{(2B + \lambda)(2C + \lambda)}{B_1C_1}, \\
q_0^2 &= \frac{4CA\lambda_0 - (C + A)^2\lambda + \lambda^2}{C_1A_1} = \frac{(2C + \lambda)(2A + \lambda)}{C_1A_1}, \\
r_0^2 &= \frac{4AB\lambda_0 - (A + B)^2\lambda + \lambda^2}{A_1B_1} = \frac{(2A + \lambda)(2B + \lambda)}{A_1B_1}.
\end{aligned}
\end{equation}
If one adds the first three of equations (3) after having multiplied them respectively by $x_0$, $y_0$, $z_0$, one finds
\begin{equation}
x_0(Ap_0 + A_1q_0r_0) + y_0(Bq_0 + B_1r_0p_0) + z_0(Cr_0 + C_1p_0q_0) = 0.
\end{equation}
This equation serves to determine the quantity $\lambda$.
The three quantities $f_0$, $g_0$, $h_0$ must satisfy the three following equations:
\begin{align*}
p_0f_0 + q_0g_0 + r_0h_0 &= 0, \\
Ap_0f_0 + Bq_0g_0 + Cr_0h_0 &= 0, \\
x_0f_0 + y_0g_0 + z_0h_0 &= \lambda,
\end{align*}
from which one obtains, by setting
\[
\mu = -(A_1x_0q_0r_0 + B_1y_0r_0p_0 + C_1z_0p_0q_0) = Ax_0p_0 + By_0q_0 + Cz_0r_0,
\]

\begin{equation}
\begin{aligned}
f_0 &= -A_1q_0r_0\frac{\lambda}{\mu}, \\
g_0 &= -B_1r_0p_0\frac{\lambda}{\mu}, \\
h_0 &= -C_1p_0q_0\frac{\lambda}{\mu}.
\end{aligned}
\end{equation}
Equations (6) and (7) determine the quantities $p_0$, $q_0$, $r_0$ up to sign.
If one sets
\[
a = \sqrt{\frac{2A + \lambda}{A_1}}, \quad b = \sqrt{\frac{2B + \lambda}{B_1}}, \quad c = \sqrt{\frac{2C + \lambda}{C_1}}
\]
fixing the sign of the quantities $a$, $b$, $c$ in an arbitrary manner, one finds that to satisfy all equations in (3) one must set
\begin{align*}
p_0 &= bc, \\
q_0 &= ca, \\
r_0 &= -ab.
\end{align*}
Equation (7) can therefore be written in the following manner
\[
x_0(A + \lambda)p_0 + y_0(B + \lambda)q_0 + z_0(C + \lambda)r_0 = 0
\]
or
\[
x_0(A + \lambda)bc + y_0(B + \lambda)ca - z_0(C + \lambda)ab = 0.
\]

\textbf{Second case.} Suppose now that one of the three quantities $A_1$, $B_1$, $C_1$ (for example $C_1$) is equal to zero. (In this case one can also always suppose $y_0 = 0$.) To satisfy the two equations
\begin{align*}
A(p_0f_0 + q_0g_0) + Cr_0h_0 &= 0, \\
p_0f_0 + q_0g_0 + r_0h_0 &= 0,
\end{align*}
it is necessary that
\[
r_0h_0 = 0.
\]
Equations (3) therefore admit two systems of solutions, which can be written, by denoting by $\varepsilon \pm 1$:
\begin{align*}
\text{I.} \quad p_0 &= \varepsilon i \cdot \frac{2C}{A - 2C} \cdot \frac{z_0}{x_0}, & f_0 &= -\frac{2C}{x_0}, \\
q_0 &= \varepsilon ip_0, & g_0 &= -i\varepsilon\frac{2C}{x_0}, \\
r_0 &= 2\varepsilon i, & h_0 &= 0.
\end{align*}

\begin{align*}
\text{II.} \quad p_0 &= 0, & f_0 &= -\frac{2A}{x_0 - i\varepsilon z_0}, \\
q_0 &= 2\varepsilon i, & g_0 &= 0, \\
r_0 &= 0, & h_0 &= \varepsilon i\frac{2A}{x_0 - i\varepsilon z_0}.
\end{align*}
Once the coefficients $p_0$, $q_0$, $r_0$, $f_0$, $g_0$, $h_0$ are determined, one obtains the coefficients
$p_m$, $q_m$, $r_m$, $f_m$, $g_m$, $h_m$ by solving the following system of linear equations:

\begin{align*}
(m - 1)Ap_m - A_1(q_0r_m + r_0q_m) - z_0g_m + y_0h_m &= P_m, \\
(m - 1)Bq_m - B_1(r_0p_m + p_0r_m) - x_0h_m + z_0f_m &= Q_m, \\
(m - 1)Cr_m - C_1(p_0q_m + q_0p_m) - y_0f_m + x_0g_m &= R_m, \\
(m - 2)f_m - r_0g_m + q_0h_m - g_0r_m + h_0q_m &= F_m, \\
(m - 2)g_m - p_0h_m + r_0f_m - h_0p_m + f_0r_m &= G_m, \\
(m - 2)h_m - q_0f_m + p_0g_m - f_0q_m + g_0p_m &= H_m.
\end{align*}
The right-hand sides of these equations are entire functions of the coefficients $p_\mu$, $q_\mu$, $r_\mu$, $f_\mu$, $g_\mu$, $h_\mu$ in which the index $\mu < m$.
In order that the series
\begin{align*}
p &= t^{-1}\sum p_m t^m, & \gamma &= t^{-2}\sum f_m t^m, \\
q &= t^{-1}\sum q_m t^m, & \gamma' &= t^{-2}\sum g_m t^m, \\
r &= t^{-1}\sum r_m t^m, & \gamma'' &= t^{-2}\sum h_m t^m,
\end{align*}
contain the sufficient number of arbitrary constants, it is necessary that the determinant of these linear equations, which is an entire function of the 6th degree in $m$, vanishes for five different positive integer values of $m$. Moreover, it is necessary that the quantities $P_m$, $Q_m$, $R_m$, $F_m$, $G_m$, $H_m$ be such that the preceding equations are compatible with one another.

In carrying out the calculations, I have ascertained that these conditions are not satisfied in the general case; but that, besides the two already known cases, they are still satisfied in a new case, where the constants satisfy the following equations:
\[
A = B = 2C, \quad z_0 = 0.
\]
It is this case that I propose to study in the following paragraphs.

\section*{§ 2.}
\setcounter{equation}{0}

In the case that we are going to consider, one has
\[
A = B = 2C, \quad z_0 = 0.
\]
By a rotation of the coordinate axes in the $xy$ plane and by a suitable choice of the unit of length, one can always arrange things in this case such that
\[
y_0 = 0, \quad C = 1.
\]
If we then set
\[
c_0 = Mgx_0,
\]
the differential equations that we have to consider are the following:

\begin{equation}
\begin{aligned}
2\frac{dp}{dt} &= qr, & \frac{d\gamma}{dt} &= r\gamma' - q\gamma'', \\
2\frac{dq}{dt} &= -pr - c_0\gamma'', & \frac{d\gamma'}{dt} &= p\gamma'' - r\gamma, \\
\frac{dr}{dt} &= c_0\gamma', & \frac{d\gamma''}{dt} &= q\gamma - p\gamma'.
\end{aligned}
\end{equation}
The three algebraic integrals of the general case are in this case:

\begin{equation}
\begin{aligned}
2(p^2 + q^2) + r^2 &= 2c_0\gamma + 6l_1, \\
2(p\gamma + q\gamma') + r\gamma'' &= 2l, \\
\gamma^2 + \gamma'^2 + \gamma''^2 &= 1.
\end{aligned}
\end{equation}
Besides these three algebraic integrals one easily finds yet a fourth.
One has in effect (setting $i = \sqrt{-1}$),
\begin{align*}
2\frac{d(p + qi)}{dt} &= -ri(p + qi) - c_0i\gamma'', \\
\frac{d(\gamma + \gamma'i)}{dt} &= -ri(\gamma + \gamma'i) + \gamma''i(p + qi).
\end{align*}
Consequently
\[
\frac{d}{dt}\{(p + qi)^2 + c_0(\gamma + \gamma'i)\} = -ri\{(p + qi)^2 + c_0(\gamma + \gamma'i)\},
\]
and similarly
\[
\frac{d}{dt}\{(p - qi)^2 + c_0(\gamma - \gamma'i)\} = ri\{(p - qi)^2 + c_0(\gamma - \gamma'i)\},
\]
from which it follows that
\[
\frac{d}{dt}\log\{(p + qi)^2 + c_0(\gamma + i\gamma')\} + \frac{d}{dt}\log\{(p - qi)^2 + c_0(\gamma - i\gamma')\} = 0,
\]
and, on integrating,
\begin{equation}
\{(p + qi)^2 + c_0(\gamma + i\gamma')\}\{(p - qi)^2 + c_0(\gamma - i\gamma')\} = k^2,
\end{equation}
denoting by $k$ an arbitrary real constant.
Let us set
\begin{align*}
x_1 &= p + qi, & y_1 &= \gamma + \gamma'i, \\
x_2 &= p - qi, & y_2 &= \gamma - \gamma'i, \\
\xi_1 &= (p + qi)^2 + c_0(\gamma + i\gamma') = x_1^2 + c_0y_1, \\
\xi_2 &= (p - qi)^2 + c_0(\gamma - i\gamma') = x_2^2 + c_0y_2.
\end{align*}
Equation (3) is then written
\[
\xi_1\xi_2 = k^2.
\]

From equations (2) one obtains the values of $r^2$, $r\gamma''$, $\gamma''^2$ expressed in $x_1$, $x_2$, $\xi_1$, $\xi_2$:
\begin{equation}
\begin{aligned}
r^2 &= 6l_1 - (x_1 + x_2)^2 + \xi_1 + \xi_2, \\
c_0r\gamma'' &= 2lc_0 + x_1x_2(x_1 + x_2) - x_2\xi_1 - x_1\xi_2, \\
c_0^2\gamma''^2 &= c_0^2 - k^2 - x_1^2x_2^2 + x_2^2\xi_1 + x_1^2\xi_2,
\end{aligned}
\end{equation}
or, setting
\begin{align*}
\mathfrak{A} &= 6l_1 - (x_1 + x_2)^2, \\
\mathfrak{B} &= 2lc_0 + x_1x_2(x_1 + x_2), \\
\mathfrak{C} &= c_0^2 - k^2 - x_1^2x_2^2,
\end{align*}
\begin{equation}
    \tag{4'}
\begin{aligned}
r^2 &= \mathfrak{A} + \xi_1 + \xi_2, \\
c_0r\gamma'' &= \mathfrak{B} - x_2\xi_1 - x_1\xi_2, \\
c_0^2\gamma''^2 &= \mathfrak{C} + x_2^2\xi_1 + x_1^2\xi_2.
\end{aligned}
\end{equation}
The four quantities $x_1$, $x_2$, $\xi_1$, $\xi_2$ therefore satisfy the two following algebraic relations
\begin{equation}
\begin{aligned}
&\xi_1\xi_2 = k^2, \\
&(\mathfrak{A} + \xi_1 + \xi_2)(\mathfrak{C} + x_2^2\xi_1 + x_1^2\xi_2) - (\mathfrak{B} - x_2\xi_1 - x_1\xi_2)^2 = 0.
\end{aligned}
\end{equation}
Let us set
\begin{align*}
R(x_1) &= \mathfrak{A}x_1^2 + 2\mathfrak{B}x_1 + \mathfrak{C} = -x_1^4 + 6l_1x_1^2 + 4lc_0x_1 + c_0^2 - k^2, \\
R(x_2) &= \mathfrak{A}x_2^2 + 2\mathfrak{B}x_2 + \mathfrak{C} = -x_2^4 + 6l_1x_2^2 + 4lc_0x_2 + c_0^2 - k^2, \\
R_1(x_1x_2) &= \mathfrak{A}\mathfrak{C} - \mathfrak{B}^2 = -6l_1x_1^2x_2^2 - (c_0^2 - k^2)(x_1 + x_2)^2 \\
&\quad - 4lc_0(x_1 + x_2)x_1x_2 + 6l_1(c_0^2 - k^2) - 4l^2c_0^2.
\end{align*}
Equation (5) then becomes
\begin{equation}
    \tag{5'}
    \begin{aligned}
    &\xi_1\xi_2 = k^2, \\
    &R(x_2)\xi_1 + R(x_1)\xi_2 + R_1(x_1x_2) + k^2(x_1 - x_2)^2 = 0.
    \end{aligned}
\end{equation}
Let us set moreover
\begin{align*}
R(x_1x_2) &= \mathfrak{A}x_1x_2 + \mathfrak{B}(x_1 + x_2) + \mathfrak{C} \\
&= -x_1^2x_2^2 + 6l_1x_1x_2 + 2lc_0(x_1 + x_2) + c_0^2 - k^2;
\end{align*}
such that one then has the identity
\begin{equation}
R(x_1)R(x_2) - R(x_1x_2)^2 = (x_1 - x_2)^2R_1(x_1x_2).
\end{equation}
Let us set:
\begin{equation}
\begin{aligned}
W^2 &= \{R_1(x_1x_2) + k^2(x_1 - x_2)^2\}^2 - 4k^2R(x_1)R(x_2) \\
&= \{R_1(x_1x_2) + k^2(x_1 - x_2)^2\}^2 - 4k^2\{R(x_1x_2)^2 + (x_1 - x_2)^2R_1(x_1x_2)\} \\
&= \{R_1(x_1x_2) - k^2(x_1 - x_2)^2\}^2 - 4k^2R(x_1x_2)^2 \\
&= \{R_1(x_1x_2) - k^2(x_1 - x_2)^2 - 2kR(x_1x_2)\}\{R_1(x_1x_2) - k^2(x_1 - x_2)^2 + 2kR(x_1x_2)\}.
\end{aligned}
\end{equation}
Since by virtue of identity (6) one has
\begin{align*}
k^2(x_1 - x_2)^2 + 2kR(x_1x_2) - R_1(x_1x_2) ={}& \frac{1}{(x_1 - x_2)^2}\{k^2(x_1 - x_2)^4 + 2k(x_1 - x_2)^2R(x_1x_2) + \\
&R(x_1x_2)^2 - R(x_1)R(x_2)\} \\
={}& \frac{1}{(x_1 - x_2)^2}\{[k(x_1 - x_2)^2 + R(x_1x_2)]^2 - R(x_1)R(x_2)\} \\
={}& (x_1 - x_2)^2\left\{\frac{R(x_1x_2) - \sqrt{R(x_1)}\sqrt{R(x_2)}}{(x_1 - x_2)^2} + k\right\} \\
& \cdot\left\{\frac{R(x_1x_2) + \sqrt{R(x_1)}\sqrt{R(x_2)}}{(x_1 - x_2)^2} + k\right\},
\end{align*}
and similarly
\begin{align*}
k^2(x_1 - x_2)^2 - 2kR(x_1x_2) - R_1(x_1x_2) =& (x_1 - x_2)^2\left\{\frac{R(x_1x_2) - \sqrt{R(x_1)}\sqrt{R(x_2)}}{(x_1 - x_2)^2} - k\right\} \\
&\cdot \left\{\frac{R(x_1x_2) + \sqrt{R(x_1)}\sqrt{R(x_2)}}{(x_1 - x_2)^2} - k\right\},
\end{align*}
one can write
\begin{equation}
    \tag{7'}
\begin{aligned}
W^2 ={}& (x_1 - x_2)^4\left\{\frac{R(x_1x_2) - \sqrt{R(x_1)}\sqrt{R(x_2)}}{(x_1 - x_2)^2} - k\right\}\left\{\frac{R(x_1x_2) - \sqrt{R(x_1)}\sqrt{R(x_2)}}{(x_1 - x_2)^2} + k\right\} \\
& \times \left\{\frac{R(x_1x_2) + \sqrt{R(x_1)}\sqrt{R(x_2)}}{(x_1 - x_2)^2} + k\right\}\left\{\frac{R(x_1x_2) + \sqrt{R(x_1)}\sqrt{R(x_2)}}{(x_1 - x_2)^2} - k\right\},
\end{aligned}
\end{equation}
and one then has
\begin{equation}
\begin{aligned}
\xi_1 &= -\frac{R_1(x_1x_2) + k^2(x_1 - x_2)^2 - W}{2R(x_2)}, \\
\xi_2 &= -\frac{R_1(x_1x_2) + k^2(x_1 - x_2)^2 + W}{2R(x_1)}.
\end{aligned}
\end{equation}

Instead of the two variables $x_1$ and $x_2$ let us now introduce two new variables $s_1$ and $s_2$, defined by the equations
\begin{equation}
\begin{aligned}
s_1 &= \frac{R(x_1x_2) - \sqrt{R(x_1)R(x_2)}}{2(x_1 - x_2)^2} + \frac{l_1}{2}, \\
s_2 &= \frac{R(x_1x_2) + \sqrt{R(x_1)R(x_2)}}{2(x_1 - x_2)^2} + \frac{l_1}{2};
\end{aligned}
\end{equation}
where $s_1$ and $s_2$ are, as one sees, the two roots of the algebraic equation of the second degree
\[
(x_1 - x_2)^2\left(s - \frac{1}{2}l_1\right)^2 - R(x_1x_2)\left(s - \frac{1}{2}l_1\right) - \frac{1}{4}R_1(x_1x_2) = 0.
\]
One then has
\[
W^2 = (x_1 - x_2)^4(2s_1 - l_1 - k)(2s_1 - l_1 + k)(2s_2 - l_1 - k)(2s_2 - l_1 + k);
\]
or, setting
\[
k_1 = \frac{l_1 + k}{2}, \quad k_2 = \frac{l_1 - k}{2},
\]
one has
\begin{equation}
W^2 = 16(x_1 - x_2)^4(s_1 - k_1)(s_2 - k_1)(s_1 - k_2)(s_2 - k_2),
\end{equation}
\begin{equation}
\begin{cases}
\displaystyle \xi_1 = \frac{(x_1 - x_2)^2}{R(x_2)}\{[\sqrt{(s_1 - k_1)(s_2 - k_1)} + \sqrt{(s_1 - k_2)(s_2 - k_2)}]^2 - k^2\}, \\[2ex]
\displaystyle \xi_2 = \frac{(x_1 - x_2)^2}{R(x_1)}\{[\sqrt{(s_1 - k_1)(s_2 - k_1)} - \sqrt{(s_1 - k_2)(s_2 - k_2)}]^2 - k^2\}.
\end{cases}
\end{equation}

Let us now establish the differential equations which define these two new variables $s_1$ and $s_2$ as functions of time.
Setting
\begin{equation}
\begin{aligned}
g_2 &= k^2 - c_0^2 + 3l_1^2, \\
g_3 &= l_1(k^2 - c_0^2 - l_1^2) + l^2c_0^2, \\
S_1 &= 4s_1^3 - g_2s_1 - g_3, \\
S_2 &= 4s_2^3 - g_2s_2 - g_3,
\end{aligned}
\end{equation}
one finds
\begin{align*}
\sqrt{S_1} &= \frac{R(x_1) - \frac{1}{4}(x_1 - x_2)R'(x_1)}{(x_1 - x_2)^3}\sqrt{R(x_2)} - \frac{R(x_2) + \frac{1}{4}(x_1 - x_2)R'(x_2)}{(x_1 - x_2)^3}\sqrt{R(x_1)}, \\
\sqrt{S_2} &= \frac{R(x_1) - \frac{1}{4}(x_1 - x_2)R'(x_1)}{(x_1 - x_2)^3}\sqrt{R(x_2)} + \frac{R(x_2) + \frac{1}{4}(x_1 - x_2)R'(x_2)}{(x_1 - x_2)^3}\sqrt{R(x_1)},
\end{align*}

\begin{equation}
\begin{cases}
\displaystyle \frac{ds_1}{\sqrt{S_1}} = \frac{dx_1}{\sqrt{R(x_1)}} + \frac{dx_2}{\sqrt{R(x_2)}}, \\[2ex]
\displaystyle \frac{ds_2}{\sqrt{S_2}} = -\frac{dx_1}{\sqrt{R(x_1)}} + \frac{dx_2}{\sqrt{R(x_2)}}.
\end{cases}
\end{equation}
From the differential equations in (1) it follows, since $x_1 = p + qi$ and $x_2 = p - qi$, that

\begin{equation}
\begin{aligned}
2i\frac{dx_1}{dt} &= rx_1 + c_0\gamma'', \\
-2i\frac{dx_2}{dt} &= rx_2 + c_0\gamma''.
\end{aligned}
\end{equation}
Consequently
\[
-4\left(\frac{dx_1}{dt}\right)^2 = r^2x_1^2 + 2rc_0\gamma''x_1 + c_0^2\gamma''^2
\]
or, writing in place of $r^2$, $c_0r\gamma''$, $c_0^2\gamma''^2$ their values, defined by equations (4'),
\begin{align*}
-4\left(\frac{dx_1}{dt}\right)^2 &= x_1^2(\mathfrak{A} + \xi_1 + \xi_2) + 2x_1(\mathfrak{B} - x_2\xi_1 - x_1\xi_2) + \mathfrak{C} + x_2^2\xi_1 + x_1^2\xi_2 \\
&= R(x_1) + (x_1 - x_2)^2\xi_1.
\end{align*}
In the same manner one finds
\[
-4\left(\frac{dx_2}{dt}\right)^2 = R(x_2) + (x_1 - x_2)^2\xi_2
\]
and
\begin{align*}
4\frac{dx_1}{dt}\cdot\frac{dx_2}{dt} &= r^2x_1x_2 + rc_0\gamma''(x_1 + x_2) + c_0^2\gamma''^2 \\
&= x_1x_2(\mathfrak{A} + \xi_1 + \xi_2) + (x_1 + x_2)(\mathfrak{B} - x_2\xi_1 - x_1\xi_2) + \mathfrak{C} + x_2^2\xi_1 + x_1^2\xi_2 \\
&= R(x_1x_2).
\end{align*}
Squaring both sides of the first of equations (13), one finds
\begin{align*}
-4\cdot\frac{1}{S_1}\left(\frac{ds_1}{dt}\right)^2 &= -\frac{4}{R(x_1)}\left(\frac{dx_1}{dt}\right)^2 - \frac{4}{R(x_2)}\left(\frac{dx_2}{dt}\right)^2 - \frac{8}{\sqrt{R(x_1)}\sqrt{R(x_2)}}\frac{dx_1}{dt}\frac{dx_2}{dt} \\
&= 2 + (x_1 - x_2)^2\left(\frac{\xi_1}{R(x_1)} + \frac{\xi_2}{R(x_2)}\right) - 2\frac{R(x_1x_2)}{\sqrt{R(x_1)}\sqrt{R(x_2)}},
\end{align*}
or, by virtue of equations (5),
\begin{align*}
-4\left(\frac{1}{\sqrt{S_1}}\frac{ds_1}{dt}\right)^2 =& 2 - \frac{(x_1 - x_2)^2}{R(x_1)\cdot R(x_2)}\{R_1(x_1x_2) + k^2(x_1 - x_2)^2\} - \frac{2R(x_1x_2)}{\sqrt{R(x_1)}\cdot\sqrt{R(x_2)}} \\
=& \frac{R(x_1)\cdot R(x_2) + R(x_1x_2)^2 - 2R(x_1x_2)\sqrt{R(x_1)}\sqrt{R(x_2)} - k^2(x_1 - x_2)^4}{R(x_1)R(x_2)} \\
=& \frac{(x_1 - x_2)^4}{R(x_1)\cdot R(x_2)}\left\{\frac{\left[R(x_1x_2) - \sqrt{R(x_1)}\sqrt{R(x_2)}\right]^2}{(x_1 - x_2)^4} - k^2\right\} \\
=& 4\frac{(x_1 - x_2)^4}{R(x_1)R(x_2)}\left(\frac{R(x_1x_2) - \sqrt{R(x_1)}\sqrt{R(x_2)}}{2(x_1 - x_2)^2} - \frac{1}{2}k\right) \\
&\cdot\left(\frac{R(x_1x_2) - \sqrt{R(x_1)}\sqrt{R(x_2)}}{2(x_1 - x_2)^2} + \frac{1}{2}k\right) \\
=& 4\frac{(x_1 - x_2)^4}{R(x_1)\cdot R(x_2)}(s_1 - k_1)(s_1 - k_2).
\end{align*}
But one has
\[
s_2 - s_1 = \frac{\sqrt{R(x_1)}\sqrt{R(x_2)}}{(x_1 - x_2)^2}.
\]
Consequently
\[
-\frac{1}{S_1}\left(\frac{ds_1}{dt}\right)^2 = \frac{(s_1 - k_1)(s_1 - k_2)}{(s_1 - s_2)^2}.
\]
Let us set
\[
R_1(s) = -S(s - k_1)(s - k_2) = -4(s - e_1)(s - e_2)(s - e_3)(s - k_1)(s - k_2);
\]
one then has
\[
\frac{ds_1}{\sqrt{R_1(s_1)}} = \frac{dt}{s_1 - s_2}.
\]
One finds in entirely the same manner
\[
\frac{ds_2}{\sqrt{R_1(s_2)}} = \frac{dt}{s_2 - s_1}.
\]
From which it follows that
\begin{equation}
\begin{aligned}
0 &= \frac{ds_1}{\sqrt{R_1(s_1)}} + \frac{ds_2}{\sqrt{R_1(s_2)}}, \\
dt &= \frac{s_1ds_1}{\sqrt{R_1(s_1)}} + \frac{s_2ds_2}{\sqrt{R_1(s_2)}}.
\end{aligned}
\end{equation}
With $R_1(s)$ a polynomial of the fifth degree, and the roots of the equation $R_1(s) = 0$ all different from each other, the differential equations in (15) lead us to hyperelliptic functions, or in other words, to the functions of Mr. ROSENHAIN.

Let us first seek the expressions of the six quantities $p$, $q$, $r$, $\gamma$, $\gamma'$, $\gamma''$ expressed with the aid of the two quantities $s_1$ and $s_2$. 
We shall arrive at this with the aid of formulas, which I borrow from an unpublished course of Mr. WEIERSTRASS on elliptic functions, and of which a part are also developed in the \textit{Traité des fonctions elliptiques} etc. of Mr. HALPHEN and in the \textit{Formeln und Lehrsätze} etc. of Mr. SCHWARTZ.

\section*{§ 3.}
\setcounter{equation}{0}

Let
\[
R(x) = Ax^4 + 4Bx^3 + 6Cx^2 + 4B'x + A'
\]
be a polynomial of the fourth degree in the variable $x$. The coefficients $A$, $B$, $C$, $B'$, $A'$ are constants, subject to the condition that $R(x)$ has no quadratic divisor. Let $u$ be a second variable, linked to $x$ by the differential equation
\begin{equation}
{du} = \frac{dx}{\sqrt{R(x)}}.
\end{equation}
The most general relation between $u$ and $x$ can be expressed in the following manner.
Let us set
\begin{equation}
\begin{aligned}
g_2 &= AA' - 4BB' + 3C^2, \\
g_3 &= ACA' + 2BCB' - AB'^2 - A'B^2 - C^3, \\
D &= B^2 - AC, \\
E &= A^2B' - 3ABC + 2B^3;
\end{aligned}
\end{equation}
so we have
\begin{equation}
4\left(\frac{D}{A}\right)^3 - g_2\frac{D}{A} - g_3 = \frac{E^2}{A^2}.
\end{equation}
We shall denote by $\wp(u)$ the function $\wp(u, g_2, g_3)$, and by $\bar{\wp}(u)$ the function $\wp(u, g_2, -g_3)$\footnote{Voir HALPHEN, \textit{Traité des fonctions elliptiques}. SCHWARTZ, \textit{Formeln und Lehrsätze zum Gebrauche der elliptischen Functionen}.}. By virtue of equation (3) one can define an argument $w$ such that, the sign of $\sqrt{A}$ being fixed arbitrarily, the two equations
\begin{equation}
\wp(w) = \frac{D}{A}, \quad \wp'(w) = -\frac{E}{A\sqrt{A}}
\end{equation}
hold simultaneously. We set
\begin{equation}
\begin{aligned}
\varphi_0(u) &= -\frac{B}{A} + \frac{1}{2\sqrt{A}}\frac{\wp u + \wp'w}{\wp u - \wp w} \\
&= -\frac{B}{A} + \frac{1}{\sqrt{A}}\left\{\frac{\sigma'(u - w)}{\sigma(u - w)} - \frac{\sigma'(u)}{\sigma(u)} + \frac{\sigma'(w)}{\sigma(w)}\right\}.
\end{aligned}
\end{equation}
The most general expression of $x$ as a function of $u$ is then
\begin{equation}
x = \varphi_0(u - u_0).
\end{equation}
where $u_0$ denotes an arbitrary constant. (For $u = u_0$, $x = \infty$ and $\frac{1}{x^2}\frac{dx}{du} = \sqrt{A}$.) One has
\begin{equation}
\varphi_0(-u + w) = \varphi_0(u), \quad \varphi_0\left(u + \frac{w}{2}\right) = \varphi_0\left(-u + \frac{w}{2}\right)
\end{equation}
and consequently, if we set
\begin{equation}
\begin{aligned}
\varphi(u, w) = \varphi_0\left(u + \frac{w}{2}\right) &= -\frac{B}{A} + \frac{1}{\sqrt{A}}\left\{\frac{\sigma'\left(u - \frac{w}{2}\right)}{\sigma\left(u - \frac{w}{2}\right)} - \frac{\sigma'\left(u + \frac{w}{2}\right)}{\sigma\left(u + \frac{w}{2}\right)} + \frac{\sigma'(w)}{\sigma(w)}\right\} \\
&= -\frac{B}{A} + \frac{1}{\sqrt{A}}\left\{\frac{\wp'\left(\frac{w}{2}\right)}{\wp u - \wp\left(\frac{w}{2}\right)} + \frac{1}{2}\frac{\wp''\left(\frac{w}{2}\right)}{\wp'\left(\frac{w}{2}\right)}\right\}.
\end{aligned}
\end{equation}
then $\varphi(u, w)$ is an even function of $u$.
Equations (4) are unchanged if one adds to $w$ any period of the function $\wp(u)$. 
Consequently, if we denote by $(2\tilde{\omega}, 2\tilde{\omega}')$ any pair of primitive periods of the function $\wp(u)$, and if of all the values of $w$ (infinite in number) which satisfy equations (4), we fix one arbitrarily, then the four expressions
\[
\varphi(u, w), \quad \varphi(u, w + 2\tilde{\omega}), \quad \varphi(u, w + 2\tilde{\omega} + 2\tilde{\omega}'), \quad \varphi(u, w + 2\tilde{\omega}')
\]
will represent four different even functions of $u$, which, when substituted for $x$, all satisfy the differential equation (1). 
In effect, if $2\Omega$ and $2\Omega'$ are two arbitrary periods of the function $\bar{\wp}(u)$, it follows from equation (8) that
\begin{equation}
\begin{aligned}
\varphi(u, w + 2\Omega) &= \varphi(u + \Omega, w), \\
\varphi(u, w + 2\Omega') &= \varphi(u + \Omega', w).
\end{aligned}
\end{equation}
For the equation
\[
\varphi(u, w + 2\Omega) = \varphi(u, w + 2\Omega')
\]
to hold for all values of $u$, it is therefore necessary that one has
\begin{align*}
\varphi(u + \Omega, w) &= \varphi(u + \Omega', w), \\
\varphi(u + \Omega' - \Omega, w) &= \varphi(u, w),
\end{align*}
that is to say $\Omega' - \Omega$ must be a period of the function $\varphi(u)$ and consequently, (in view of the last of equations (8)) $\Omega' - \Omega$ must also be a period of $\wp u$.
If this is not the case, the two functions $\varphi(u, w + 2\Omega)$ and $\varphi(u, w + 2\Omega')$ are not identical. 
One concludes therefore that the four aforementioned functions,
\[
\varphi(u, w), \quad \varphi(u, w + 2\tilde{\omega}), \quad \varphi(u, w + 2\tilde{\omega} + 2\tilde{\omega}'), \quad \varphi(u, w + 2\tilde{\omega}'),
\]
are all distinct.
But there exist in total only four even functions of $u$, which, when substituted for $x$, satisfy the differential equation (1); these four functions are distinguished by the fact that for $u = 0$ each of them becomes equal to one of the four roots of the equation
\[
R(x) = 0.
\]
These four functions must therefore be identical with the four aforementioned even functions, and, by virtue of equations (9) one can also write them in the following manner:
\[
\varphi(u, w), \quad \varphi(u + \tilde{\omega}, w), \quad \varphi(u + \tilde{\omega} + \tilde{\omega}', w), \quad \varphi(u + \tilde{\omega}', w).
\]
If we set
\begin{align*}
a &= \varphi(0, w), & a_1 &= \varphi(\tilde{\omega}, w) = \varphi(0, w + 2\tilde{\omega}), \\
a_2 &= \varphi(\tilde{\omega} + \tilde{\omega}', w) = \varphi(0, w + 2\tilde{\omega} + 2\tilde{\omega}'), & a_3 &= \varphi(\tilde{\omega}', w) = \varphi(0, w + 2\tilde{\omega}'),
\end{align*}
the four quantities $a$, $a_1$, $a_2$, $a_3$ represent the 4 different roots of the equation $R(x) = 0$ and one has
\begin{align*}
a &= -\frac{B}{A} + \frac{1}{2\sqrt{A}}\frac{\wp''\left(\frac{w}{2}\right)}{\wp'\left(\frac{w}{2}\right)}, &
a_2 &= -\frac{B}{A} + \frac{1}{2\sqrt{A}}\frac{\wp''\left(\frac{w}{2} + \tilde{\omega} + \tilde{\omega}'\right)}{\wp'\left(\frac{w}{2} + \tilde{\omega} + \tilde{\omega}'\right)}, \\
a_1 &= -\frac{B}{A} + \frac{1}{2\sqrt{A}}\frac{\wp''\left(\frac{w}{2} + \tilde{\omega}\right)}{\wp'\left(\frac{w}{2} + \tilde{\omega}\right)}, &
a_3 &= -\frac{B}{A} + \frac{1}{2\sqrt{A}}\frac{\wp''\left(\frac{w}{2} + \tilde{\omega}'\right)}{\wp'\left(\frac{w}{2} + \tilde{\omega}'\right)}.
\end{align*}
Setting $\wp(\tilde{\omega}) = e_1$, $\wp(\tilde{\omega} + \tilde{\omega}') = e_2$, $\wp(\tilde{\omega}') = e_3$, the following relations hold for each value of the argument $u$,
\begin{align*}
\wp\left(\frac{u}{2}\right) &= -\frac{(e_2^2 - e_3^2)\sigma_1 u + (e_3^2 - e_1^2)\sigma_2 u + (e_1^2 - e_2^2)\sigma_3 u}{(e_2 - e_3)\sigma_1 u + (e_3 - e_1)\sigma_2 u + (e_1 - e_2)\sigma_3 u}, \\
\wp'\left(\frac{u}{2}\right) &= \frac{-2(e_2 - e_3)(e_3 - e_1)(e_1 - e_2)\sigma u}{(e_2 - e_3)\sigma_1 u + (e_3 - e_1)\sigma_2 u + (e_1 - e_2)\sigma_3 u}, \\
\frac{\wp''\left(\frac{u}{2}\right)}{\wp'\left(\frac{u}{2}\right)} &= -2\cdot\frac{(e_2 - e_3)\sigma_2 u\sigma_3 u + (e_3 - e_1)\sigma_3 u\sigma_1 u + (e_1 - e_2)\sigma_1 u\sigma_2 u}{\sigma u[(e_2 - e_3)\sigma_1 u + (e_3 - e_1)\sigma_2 u + (e_1 - e_2)\sigma_3 u]} \\
&= -2\frac{\sigma_1 u + \sigma_2 u + \sigma_3 u}{\sigma u}.\footnotemark
\end{align*}
\footnotetext{See Schwarz, \textit{Formeln und Lehrs\"atze}.}
One has moreover
\[
\frac{\sigma_1(w)}{\sigma(w)} = \sqrt{\frac{D}{A} - e_1}, \quad \frac{\sigma_2(w)}{\sigma(w)} = \sqrt{\frac{D}{A} - e_2}, \quad \frac{\sigma_3(w)}{\sigma(w)} = \sqrt{\frac{D}{A} - e_3}.
\]
For each determined value of $w$, the corresponding signs of the square roots in these equations are also perfectly determined, and one has
\[
\sqrt{\frac{D}{A} - e_1} \cdot \sqrt{\frac{D}{A} - e_2} \cdot \sqrt{\frac{D}{A} - e_3} = -\frac{1}{2}\wp'(w) = \frac{E}{2A\sqrt{A}}.
\]
By virtue of these equations, one therefore has
\begin{align*}
-\wp\left(\frac{w}{2}\right) &= \frac{(e_2^2 - e_3^2)\sqrt{\frac{D}{A} - e_1} + (e_3^2 - e_1^2)\sqrt{\frac{D}{A} - e_2} + (e_1^2 - e_2^2)\sqrt{\frac{D}{A} - e_3}}{(e_2 - e_3)\sqrt{\frac{D}{A} - e_1} + (e_3 - e_1)\sqrt{\frac{D}{A} - e_2} + (e_1 - e_2)\sqrt{\frac{D}{A} - e_3}}, \\
\wp'\left(\frac{w}{2}\right) &= \frac{-2(e_2 - e_3)(e_3 - e_1)(e_1 - e_2)}{(e_2 - e_3)\sqrt{\frac{D}{A} - e_1} + (e_3 - e_1)\sqrt{\frac{D}{A} - e_2} + (e_1 - e_2)\sqrt{\frac{D}{A} - e_3}}, \\
\frac{1}{2\sqrt{A}}\frac{\wp''\left(\frac{w}{2}\right)}{\wp'\left(\frac{w}{2}\right)} &= -\frac{1}{\sqrt{A}}\left\{\sqrt{\frac{D}{A} - e_1} + \sqrt{\frac{D}{A} - e_2} + \sqrt{\frac{D}{A} - e_3}\right\}.
\end{align*}
Consequently the root $a$ of the equation $R(x) = 0$, corresponding to the particular value of $w$ that we have chosen, is given by the formula:
\begin{align*}
a &= -\frac{1}{\sqrt{A}}\frac{\sigma_1(w) + \sigma_2(w) + \sigma_3(w)}{\sigma(w)} - \frac{B}{A} \\
&= -\frac{B}{A} - \frac{1}{\sqrt{A}}\left\{\sqrt{\frac{D}{A} - e_1} + \sqrt{\frac{D}{A} - e_2} + \sqrt{\frac{D}{A} - e_3}\right\}.
\end{align*}
If one sets
\begin{align*}
h_0 &= \frac{\sigma_1(w) + \sigma_2(w) + \sigma_3(w)}{\sigma(w)} \\
&= \sqrt{\frac{D}{A} - e_1} + \sqrt{\frac{D}{A} - e_2} + \sqrt{\frac{D}{A} - e_3}, \\
h_1 &= \frac{(e_2 - e_3)\sigma_1(w) + (e_3 - e_1)\sigma_2(w) + (e_1 - e_2)\sigma_3(w)}{\sigma(w)} \\
&= (e_2 - e_3)\sqrt{\frac{D}{A} - e_1} + (e_3 - e_1)\sqrt{\frac{D}{A} - e_2} + (e_1 - e_2)\sqrt{\frac{D}{A} - e_3}, \\
h_2 &= \frac{(e_2^2 - e_3^2)\sigma_1(w) + (e_3^2 - e_1^2)\sigma_2(w) + (e_1^2 - e_2^2)\sigma_3(w)}{\sigma(w)} \\
&= (e_2^2 - e_3^2)\sqrt{\frac{D}{A} - e_1} + (e_3^2 - e_1^2)\sqrt{\frac{D}{A} - e_2} + (e_1^2 - e_2^2)\sqrt{\frac{D}{A} - e_3},
\end{align*}
one therefore has
\begin{align*}
a &= -\frac{B}{A} - \frac{h_0}{\sqrt{A}}, \\
\varphi(u, w) &= -\frac{B}{A} - \frac{h_0}{\sqrt{A}} - \frac{2}{\sqrt{A}}\frac{(e_2 - e_3)(e_3 - e_1)(e_1 - e_2)}{h_1\wp u + h_2}.
\end{align*}

By setting
\begin{align*}
\tilde{\omega} &= \omega_1, \\
\tilde{\omega} + \tilde{\omega}' &= \omega_2, \\
\tilde{\omega}' &= \omega_3
\end{align*}
and by designating by $\lambda$, $\mu$, $\nu$ the numbers 1, 2, 3 in any order, one has
\begin{align*}
\frac{\sigma_\lambda(u + 2\omega_\lambda)}{\sigma(u + 2\omega_\lambda)} &= \frac{\sigma_\lambda u}{\sigma u}, \\
\frac{\sigma_\lambda(u + 2\omega_\mu)}{\sigma(u + 2\omega_\mu)} &= -\frac{\sigma_\lambda u}{\sigma u}, \\
\frac{\sigma_\lambda(u + 2\omega_\nu)}{\sigma(u + 2\omega_\nu)} &= -\frac{\sigma_\lambda u}{\sigma u}.
\end{align*}
One therefore obtains, by writing in the preceding expressions $w + 2\omega_1$, $w + 2\omega_2$, $w + 2\omega_3$ in place of $w$ and by setting \hypertarget{TN1-ref}{}\hypertarget{TN2-ref}{}\tnref{TN1}, \tnref{TN2}:

\begin{equation}
    \tag{14}\label{eq:14}
\begin{aligned}
h_0' &= \frac{\sigma_1(w) - \sigma_2(w) - \sigma_3(w)}{\sigma(w)}, \\
h_1' &= \frac{(e_2 - e_3)\sigma_1(w) - (e_3 - e_1)\sigma_2(w) - (e_1 - e_2)\sigma_3(w)}{\sigma(w)}, \\
h_2' &= \frac{(e_2^2 - e_3^2)\sigma_1(w) - (e_3^2 - e_1^2)\sigma_2(w) - (e_1^2 - e_2^2)\sigma_3(w)}{\sigma(w)}, \\
h_0'' &= \frac{-\sigma_1(w) + \sigma_2(w) - \sigma_3(w)}{\sigma(w)}, \\
h_1'' &= \frac{-(e_2 - e_3)\sigma_1(w) + (e_3 - e_1)\sigma_2(w) - (e_1 - e_2)\sigma_3(w)}{\sigma(w)}, \\
h_2'' &= \frac{-(e_2^2 - e_3^2)\sigma_1(w) + (e_3^2 - e_1^2)\sigma_2(w) - (e_1^2 - e_2^2)\sigma_3(w)}{\sigma(w)}, \\
h_0''' &= \frac{-\sigma_1(w) - \sigma_2(w) + \sigma_3(w)}{\sigma(w)}, \\
h_1''' &= \frac{-(e_2 - e_3)\sigma_1(w) - (e_3 - e_1)\sigma_2(w) + (e_1 - e_2)\sigma_3(w)}{\sigma(w)}, \\
h_2''' &= \frac{-(e_2^2 - e_3^2)\sigma_1(w) - (e_3^2 - e_1^2)\sigma_2(w) + (e_1^2 - e_3^2)\sigma_3(w)}{\sigma(w)}
\end{aligned}
\end{equation}
the following
\begin{equation}
    \tag{15} 
\begin{aligned}
a &= -\frac{B}{A} - \frac{h_0}{\sqrt{A}}, \\
\varphi(u, w) &= -\frac{B}{A} - \frac{h_0}{\sqrt{A}} - \frac{2}{\sqrt{A}}\frac{(e_2 - e_3)(e_3 - e_1)(e_1 - e_2)}{h_1\wp u + h_2}, \\
a_1 &= -\frac{B}{A} - \frac{h_0'}{\sqrt{A}}, \\
\varphi(u + \omega_1, w) &= -\frac{B}{A} - \frac{h_0'}{\sqrt{A}} - \frac{2}{\sqrt{A}}\frac{(e_2 - e_3)(e_3 - e_1)(e_1 - e_2)}{h_1'\wp u + h_2'}, \\
a_2 &= -\frac{B}{A} - \frac{h_0''}{\sqrt{A}}, \\
\varphi(u + \omega_2, w) &= -\frac{B}{A} - \frac{h_0''}{\sqrt{A}} - \frac{2}{\sqrt{A}}\frac{(e_2 - e_3)(e_3 - e_1)(e_1 - e_2)}{h_1''\wp u + h_2''}, \\
a_3 &= -\frac{B}{A} - \frac{h_0'''}{\sqrt{A}}, \\
\varphi(u + \omega_3, w) &= -\frac{B}{A} - \frac{h_0'''}{\sqrt{A}} - \frac{2}{\sqrt{A}}\frac{(e_2 - e_3)(e_3 - e_1)(e_1 - e_2)}{h_1'''\wp u + h_2'''}.
\end{aligned}
\end{equation}
In this manner one obtains the four even functions of $u$, which satisfy the differential equation $\frac{dx}{du} = \sqrt{R(x)}$, as well as the corresponding roots of the equation $R(x) = 0$, that is to say the values that each of these four functions takes for $u = 0$, expressed as rational functions of the following quantities:
\[
\frac{B}{A}, \quad \frac{1}{\sqrt{A}}, \quad \frac{\sigma_1(w)}{\sigma(w)} = \sqrt{\frac{D}{A} - e_1}, \quad \frac{\sigma_2(w)}{\sigma(w)} = \sqrt{\frac{D}{A} - e_2}, \quad \frac{\sigma_3(w)}{\sigma(w)} = \sqrt{\frac{D}{A} - e_3}.
\]

Until now we have subjected the constants $A$, $B$, $C$, $B'$, $A'$ (the coefficients of $R(x)$) to no condition; let us now suppose that all these coefficients are real. 
There are always in this case real positive values of $u$ which satisfy the equation $\wp'(u) = 0$; let us denote by $\omega$ the smallest of all these values. 
Similarly the equation $\bar{\wp}'(u) = 0$ can be satisfied by real positive values of $u$ and we shall denote by $\tilde{\omega}$ the smallest of these.

It is necessary at present to distinguish two cases:

\textbf{I. The quantity $g_2^3 - 27g_3^2$ is positive.}

Let us set in this case
\[
\omega_1 = \omega, \quad \omega_2 = \omega + \tilde{\omega}i, \quad \omega_3 = \tilde{\omega}i,
\]
\[
\wp(\omega_1) = e_1, \quad \wp(\omega_2) = e_2, \quad \wp(\omega_3) = e_3.
\]
The quantities $e_1$, $e_2$, $e_3$, which represent the roots of the equation
\[
4s^3 - g_2s - g_3 = 0,
\]
are in this case all real, and one has
\[
e_1 > e_2 > e_3;
\]
and $(2\omega_1, 2\omega_3)$ is a pair of primitive periods of the function $\wp(u)$.

\textbf{II. The quantity $g_2^3 - 27g_3^2$ is negative.}

Let us set in this case
\[
\omega_1 = \frac{\omega - \tilde{\omega}i}{2}, \quad \omega_2=\omega, \quad \omega_3 = \frac{\omega + \tilde{\omega}i}{2},
\]
\[
\wp(\omega_1) = e_1, \quad \wp(\omega_2) = e_2, \quad \wp(\omega_3) = e_3.
\]
where $(2\omega_1, 2\omega_3)$ also represents in this case a pair of primitive periods of the function $\wp(u)$, but $e_1$ and $e_3$ are in this case imaginary conjugate quantities, while $e_2 = -(e_1 + e_3)$ is real. The quantity $\frac{e_1 - e_3}{i}$ is positive.
If we agree to give $\sqrt{A}$ its positive value in the case where $A$ is positive, and to denote by this root the product of $i$ by a positive quantity, in the case where $A$ is negative, one can in both cases determine in the following manner a value of $w$ satisfying the equations
\[
\wp(w) = \frac{D}{A}, \quad \wp'(w) = -\frac{E}{A\sqrt{A}}.
\]
The equations
\[
\frac{E^2}{A^3} = \wp'(w)^2 = 4\left(\frac{D}{A} - e_1\right)\left(\frac{D}{A} - e_2\right)\left(\frac{D}{A} - e_3\right)
\]
show us that in case I, if $A > 0$, $\frac{D}{A}$ must be contained either in the interval $(e_1 \ldots \infty)$, or in the interval $(e_3 \ldots e_2)$; if $A < 0$, $\frac{D}{A}$ will belong either to the interval $(-\infty \ldots e_3)$ or to the interval $(e_2 \ldots e_1)$. In case II, on the contrary, $\frac{D}{A}$ must be contained either in the interval $(e_2 \ldots \infty)$, or in the interval $(-\infty \ldots e_2)$.

Let us denote by $(a \ldots b)$ a quantity contained in the real interval $(a \ldots b)$, by $(+)$ and $(-)$ positive or negative quantities respectively and by $\epsilon$ a real quantity satisfying the condition
\[
0 < \epsilon < 1.
\]
One has in case I
\begin{align*}
\text{(1)} \quad & \wp(2\epsilon\omega) = (\infty \ldots e_1), \quad
\wp'(2\epsilon\omega) = \begin{cases}
(-), & \epsilon < \frac{1}{2}, \\
0, & \epsilon = \frac{1}{2}, \\
(+), & \epsilon > \frac{1}{2},
\end{cases} \\
\text{(2)} \quad & \wp(\omega + 2\epsilon\tilde{\omega}i) = (e_1 \ldots e_2), \quad
\wp'(\omega + 2\epsilon\tilde{\omega}i) = \begin{cases}
(+)i, & \epsilon < \frac{1}{2}, \\
0, & \epsilon = \frac{1}{2}, \\
(-)i, & \epsilon > \frac{1}{2},
\end{cases} \\
\text{(3)} \quad & \wp(2\epsilon\omega + \tilde{\omega}i) = (e_3 \ldots e_2), \quad
\wp'(2\epsilon\omega + \tilde{\omega}i) = \begin{cases}
(+), & \epsilon < \frac{1}{2}, \\
0, & \epsilon = 0, \\
(-), & \epsilon > \frac{1}{2},
\end{cases} \\
\text{(4)} \quad & \wp(2\epsilon\tilde{\omega}i) = (-\infty \ldots e_3), \quad
\wp'(2\epsilon\tilde{\omega}i) = \begin{cases}
(-)i, & \epsilon < \frac{1}{2}, \\
0, & \epsilon = \frac{1}{2}, \\
(+)i, & \epsilon > \frac{1}{2}.
\end{cases}
\end{align*}

In each of these four cases, if one lets $\epsilon$ grow in a continuous manner from 0 to $\frac{1}{2}$, 
the function $\wp(u)$ traverses the entire indicated interval, increasing continually in the 3rd and in the 4th case, decreasing in the 1st and in the 2nd until $\epsilon = \frac{1}{2}$. 
If one then varies $\epsilon$ from $\frac{1}{2}$ to 1, $\wp(u)$ traverses in each case the same interval as before, but in the opposite direction. 
Two values of $\epsilon$, at equal distance from $\frac{1}{2}$, thus correspond to the same values of $\wp(u)$, but contrary values of $\wp'(u)$.

In case II one has
\begin{align*}
\text{(5)} \quad & \wp(2\epsilon\omega) = (+\infty \ldots e_2), \quad
\wp'(2\epsilon\omega) = \begin{cases}
(-), & \epsilon < \frac{1}{2}, \\
0, & \epsilon = \frac{1}{2}, \\
(+), & \epsilon > \frac{1}{2},
\end{cases} \\
\text{(6)} \quad & \wp(2\epsilon\tilde{\omega}i) = (-\infty \ldots e_2), \quad
\wp'(2\epsilon\tilde{\omega}i) = \begin{cases}
(-)i, & \epsilon < \frac{1}{2}, \\
0, & \epsilon = \frac{1}{2}, \\
(+)i, & \epsilon > \frac{1}{2}.
\end{cases}
\end{align*}
According to these equations, one can define a quantity $w$, satisfying the equations
\[
\wp(w) = \frac{D}{A}, \quad \wp'(w) = -\frac{E}{A\sqrt{A}}
\]
in the following manner:

\textbf{I. $g_2^3 - 27g_3^2$ is positive:}

1) If $A$ is positive and $\frac{D}{A}$ belongs to the interval $(+\infty \ldots e_1)$ one can set
\[
w = 2\epsilon\omega \quad (0 < \epsilon < 1),
\]
remarking that
\[
\epsilon \lessgtr \frac{1}{2}, \quad \text{according to } E \lessgtr 0.
\]

2) If $A$ is negative and $\frac{D}{A}$ is found in the interval $(-\infty \ldots e_3)$, one can set
\[
w = 2\epsilon\tilde{\omega}i \quad (0 < \epsilon < 1),
\]
and one has
\[
\epsilon \lessgtr \frac{1}{2}, \quad \text{according to } E \lessgtr 0.
\]

3) If $A$ is positive, but $\frac{D}{A}$ is found in the interval $(e_3 \ldots e_2)$ one can set
\[
w = 2\epsilon\omega + \tilde{\omega}i \quad (0 \leq \epsilon < 1),
\]
and one has
\[
\epsilon \lessgtr \frac{1}{2}, \quad \text{according to } E \lessgtr 0.
\]

4) If $A$ is negative and $\frac{D}{A}$ belongs to the interval $(e_2 \ldots e_1)$ one can set
\[
w = \omega + 2\epsilon\tilde{\omega}i \quad (0 \leq \epsilon < 1),
\]
and one has
\[
\epsilon \lessgtr \frac{1}{2}, \quad \text{according to } E \lessgtr 0.
\]

\textbf{II. $g_2^3 - 27g_3^2$ is a negative quantity.}

5) If $A$ is positive, $\frac{D}{A}$ belongs to the interval $(e_2 \ldots \infty)$ and one can set
\[
w = 2\epsilon\omega \quad (0 \leq \epsilon < 1),
\]
$\epsilon$ being $\lessgtr \frac{1}{2}$ according to $E \lessgtr 0$.

6) If $A$ is negative, $\frac{D}{A}$ is found in the interval $(-\infty \ldots e_2)$ and one can set
\[
w = 2\epsilon\tilde{\omega}i \quad (0 < \epsilon < 1),
\]
$\epsilon$ being $\lessgtr \frac{1}{2}$ according to $E \lessgtr 0$.

Cases 1) and 2) occur when all four roots of the equation $R(x) = 0$ are real.
Cases 3) and 4) occur if all these roots are imaginary.
Finally one has cases 5) and 6) when two roots are real and the two others imaginary conjugates.

In case 1) the quantities $\frac{D}{A} - e_1$, $\frac{D}{A} - e_2$, $\frac{D}{A} - e_3$ are all real and positive; in case 2) they are real and negative; consequently, by virtue of equations (14) the quantities
\[
\frac{h_0}{\sqrt{A}}, \quad \frac{h_0'}{\sqrt{A}}, \quad \frac{h_0''}{\sqrt{A}}, \quad \frac{h_0'''}{\sqrt{A}}
\]
are in both cases real.
In case 3) the quantities $\frac{D}{A} - e_1$, $\frac{D}{A} - e_2$ are negative, while $\frac{D}{A} - e_3$ is positive; consequently
\[
\frac{1}{\sqrt{A}}\sqrt{\frac{D}{A} - e_1} \quad \text{and} \quad \frac{1}{\sqrt{A}}\sqrt{\frac{D}{A} - e_2}
\]
are imaginary quantities, while $\frac{1}{\sqrt{A}}\sqrt{\frac{D}{A} - e_3}$ is a real quantity.
The quantities $\frac{h_0}{\sqrt{A}}$, $\frac{h_0'}{\sqrt{A}}$, $\frac{h_0''}{\sqrt{A}}$, $\frac{h_0'''}{\sqrt{A}}$ are therefore complex quantities. The quantities $a$ and $a_2$, as well as $a_1$ and $a_3$ are in this case complex conjugate quantities.
In case 4) the quantities $\frac{D}{A} - e_1$, $\frac{D}{A} - e_3$ are positive, while $\frac{D}{A} - e_2$ is negative; consequently
$\frac{1}{\sqrt{A}}\sqrt{\frac{D}{A} - e_2}$ is real, but $\frac{1}{\sqrt{A}}\sqrt{\frac{D}{A} - e_1}$ and $\frac{1}{\sqrt{A}}\sqrt{\frac{D}{A} - e_3}$ are imaginary.
The quantities $\frac{h_0}{\sqrt{A}}$, $\frac{h_0'}{\sqrt{A}}$, $\frac{h_0''}{\sqrt{A}}$, $\frac{h_0'''}{\sqrt{A}}$ are complex.
The roots $a$ and $a_1$, as well as $a_2$ and $a_3$ are complex conjugate quantities.

In case 5) $\frac{D}{A} - e_2$ is positive and $\frac{D}{A} - e_1$, $\frac{D}{A} - e_3$ are complex conjugate quantities. The functions $\frac{\sigma_2u}{\sigma u}$, $\frac{\sigma_3u}{\sigma u}$ take for real values of $u$ complex conjugate values, in such a manner that if $u$ and $u'$ are complex conjugate values, the corresponding values
\[
\frac{\sigma_2u}{\sigma u}, \quad \frac{\sigma_3u'}{\sigma u'}
\]
are also complex conjugate quantities.

In case 5), as well as in case 6), the quantities
\[
\frac{1}{\sqrt{A}}\frac{\sigma_1w}{\sigma w}, \quad \frac{1}{\sqrt{A}}\frac{\sigma_3w}{\sigma w}
\]
are complex conjugate quantities, while
\[
\frac{1}{\sqrt{A}}\frac{\sigma_2w}{\sigma w}
\]
is a real quantity.

The quantities $a$ and $a_2$ are therefore real in cases 5) and 6), but $a_1$ and $a_3$ are complex conjugate.

It follows therefore from the preceding discussion that, for real values of $u$, the functions
\[
\varphi(u, w), \quad \varphi(u + \omega_1, w), \quad \varphi(u + \omega_2, w), \quad \varphi(u + \omega_3, w)
\]
are all real in cases 1) and 2) and all imaginary, in cases 3) and 4). It is to be remarked that in case 3) the first and the third as well as the second and the fourth of these functions are conjugate imaginary quantities; on the contrary, in case 4) the first and the second, as well as the 3rd and the 4th of these functions are conjugate.
In cases 5) and 6), the first and the third of these functions are real, while the second and the 4th are imaginary conjugates.

\section*{§ 4.}
\setcounter{equation}{0}

To be able to apply the formulas of the preceding section to the case which occupies us, one must begin by discussing the roots of the equation
\[
R(x) = -x^4 + 6l_1x^2 + 4c_0lx + c_0^2 - k^2 = 0.
\]
Let us set
\[
k_0 = c_0^2 - k^2, \quad l_0 = c_0l.
\]
In general, if
\[
R(x) = Ax^4 + 4Bx^3 + 6Cx^2 + 4B'x + A'
\]
and one sets
\begin{align*}
g_2 &= AA' - 4BCB' + 3C^2, \\
g_3 &= ACA' + 2BCB' - AB'^2 - A'B^2 - C^3,
\end{align*}
the condition for the reality of the roots of the equation $R(x) = 0$ can be stated in the following manner.

If $G = g_2^3 - 27g_3^2 < 0$, the equation $R(x) = 0$ has two real roots and two imaginary conjugate roots.

If $G = g_2^3 - 27g_3^2 > 0$, all four roots are real or all four imaginary, with two conjugate pairs. The case with all real roots occurs if one further has
\[
B^2 - AC > 0, \quad 12(B^2 - AC)^2 - A^2g_2 > 0.
\]

But if, when $g_2^3 - 27g_3^2$ is positive, one of the above quantities is negative, all four roots of the equation $R(x) = 0$ are imaginary.

Applying this to our case. One will have
\begin{align*}
A &= -1, \quad B = 0, \quad C = l_1, \quad B' = lc_0 = l_0, \quad A' = c_0^2 - k^2 = k_0, \\
g_2 &= -k_0 + 3l_1^2, \\
g_3 &= -l_1(k_0+l_1^2) + l_0^2,
\end{align*}
\begin{align*}
B^2 - AC &= l_1, \\
12(B^2 - AC)^2 - A^2g_2 &= k_0 + 9l_1^2, \\
g_2^3 - 27g_3^2 &= (-k_0 + 3l_1^2)^3 - 27[-l_1(k_0 + l_1^2) + l_0^2]^2 \\
&= -27\left\{l_0^4 - 2l_1(k_0 + l_1^2)l_0^2 + \frac{1}{27}k_0\left(k_0+9l_1^2\right)\right\}.
\end{align*}

By virtue of the identity
\[
l_1^2(k_0 + l_1^2)^2 = \frac{1}{27}k_0(k_0 + 9l_1^2)^2 + \left(\frac{-k_0 + 3l_1^2}{3}\right)^3
\]
the roots of the quadratic equation in $l_0^2$
\[
Q(l_0^2) = l_0^4 - 2l_1(k_0 + l_1^2)l_0^2 + \frac{1}{27}k_0(k_0 + 9l_1^2)^2
\]
can be written
\begin{align*}
l_0^{'2} &= l_1(k_0 + l_1^2) + \left(\frac{-k_0 + 3l_1^2}{3}\right)^{\frac{3}{2}}, \\
l_0^{''2} &= l_1(k_0 + l_1^2) - \left(\frac{-k_0 + 3l_1^2}{3}\right)^{\frac{3}{2}},
\end{align*}
where we shall suppose the radical $\left(\frac{-k_0 + 3l_1^2}{3}\right)^{\frac{3}{2}}$ positive, if it is real.

Applying these formulas to the discussion of the roots of the equation $R(x) = 0$, one sees that one must distinguish the following cases:

\textbf{1°.} $l_1 > 0$, $k_0 > 0$, $-k_0 + 3l_1^2 < 0$.

In this case $l_0^{'2}$ and $l_0^{''2}$ are both imaginary. Consequently $G = -27Q(l_0^2)$ is negative for all real values of $l_0$, and the equation $R(x) = 0$ has two real roots and two imaginary roots.

\textbf{2°.} $l_1 > 0$, $k_0 > 0$, $-k_0 + 3l_1^2 > 0$.

$l_0^{'2}$ and $l_0^{''2}$ are both positive and $l_0^{'2} > l_0^{''2}$. Then, if
\[
l_0^2 > l_0^{'2} \quad \text{or} \quad l_0^2 < l_0^{''2},
\]
one has $G < 0$ and the equation $R(x) = 0$ has yet two real roots and two imaginary roots. If on the contrary
\[
l_0^{'2} > l_0^2 > l_0^{''2},
\]
$G$ is positive, and, since one also has $l_1 > 0$, $k_0 + 9l_1^2 > 0$, the four roots of the equation $R(x) = 0$ are real.

\textbf{3°.} $l_1 > 0$, $k_0 < 0$, $k_0 + 9l_1^2 > 0$.

The root $l_0^{'2}$ is positive, the other $l_0^{''2}$ negative. The quantity
\[
G = -27Q(l_0^2)
\]
is negative, if $l_0^2 > l_0^{'2}$, and in this case the equation $R(x) = 0$ has four real roots. If on the contrary $l_0^2 < l_0^{'2}$, $G$ is positive, and the proposed equation has only two real roots.

\textbf{4°.} $l_1 > 0$, $k_0 < 0$, $k_0 + 9l_1^2 < 0$.

As in the preceding case $G$ is positive for $l_0^2 < l_0^{'2}$, and negative for $l_0^2 > l_0^{'2}$. If it is negative, the equation $R(x) = 0$ has in this case four imaginary roots; if it is positive it has two real roots and two imaginary roots.

\textbf{5°.} $l_1 < 0$, $k_0 > 0$.

With $l_0^{'2}$ and $l_0^{''2}$ being negative or imaginary, the quantity $G = -27Q(l_0^2)$ remains negative for every real value of $l_0$; the equation $R(x) = 0$ therefore has two real roots and two imaginary roots.

\textbf{6°.} $l_1 < 0$, $k_0 < 0$.

This case is entirely conformable to the fourth.

\textbf{Summary:}

The equation $R(x) = 0$ has four real roots in the two following cases:
\begin{align*}
\text{1)} \quad &l_1 > 0, \quad k_0 > 0, \quad -k_0 + 3l_1^2 > 0, \quad l_0^{'2} > l_0^2 > l_0^{''2}, \\
&l_0^{'2} = l_1(k_0 + l_1^2) +  \left(\frac{-k_0 + 3l_1^2}{3}\right)^{\frac{3}{2}}, \\
&l_0^{''2} = l_1(k_0 + l_1^2) -  \left(\frac{-k_0 + 3l_1^2}{3}\right)^{\frac{3}{2}}. \\
\text{2)} \quad &l_1 > 0, \quad k_0 < 0, \quad k_0 + 9l_1^2 > 0, \quad l_0^2 < l_0^{'2}.
\end{align*}

The four roots of the equation $R(x) = 0$ are imaginary in the following cases:
\begin{align*}
\text{1)} \quad &l_1 > 0, \quad k_0 > 0, \quad k_0 + 9l_1^2 < 0, \quad l_0^2 < l_0^{'2}. \\
\text{2)} \quad &l_1 < 0, \quad k_0 < 0, \quad l_0^2 < l_0^{'2}.
\end{align*}

Finally, the equation $R(x) = 0$ has two real roots in the following cases:
\begin{align*}
\text{1)} \quad &l_1 > 0, \quad k_0 > 0, \quad -k_0 + 3l_1^2 < 0, \text{ for each value of } l_0. \\
\text{2)} \quad &l_1 > 0, \quad k_0 > 0, \quad -k_0 + 3l_1^2 > 0, \quad l_0^2 > l_0^{'2} \text{ or } l_0^2 < l_0^{''2}. \\
\text{3)} \quad &l_1 > 0, \quad k_0 < 0, \quad k_0 + 9l_1^2 > 0, \quad l_0^2 < l_0^{'2}. \\
\text{4)} \quad &l_1 > 0, \quad k_0 < 0, \quad k_0 + 9l_1^2 < 0, \quad l_0^2 > l_0^{'2}. \\
\text{5)} \quad &l_1 < 0, \quad k_0 > 0, \text{ for each value of } l_0. \\
\text{6)} \quad &l_1 < 0, \quad k_0 < 0, \quad l_0^2 > l_0^{'2}.
\end{align*}

\section*{§ 5.}
\setcounter{equation}{0}

I shall examine in more detail the case where all the roots of the equation
\[
R(x) = -x^4 + 6l_1x^2 + 4lc_0x + c_0^2 - k^2 = 0
\]
are real.
Using the same notations as in § 3, one therefore has in this case
\[
\wp(w) = \frac{D}{A} = -l_1, \quad \wp'(w)^2 = 4s^3 - g_2s -g_3 = -l_0^2.
\]
where $w$ is purely imaginary.
The quantities $e_1$, $e_2$, $e_3$ are real and one has
\[
e_1 > e_2 > e_3 > -l_1.
\]
The quantities
\[
\frac{\sigma_1(w)}{\sigma(w)} = \sqrt{-(l_1 + e_1)}, \quad \frac{\sigma_2(w)}{\sigma(w)} = \sqrt{-(l_1 + e_2)}, \quad \frac{\sigma_3(w)}{\sigma(w)} = \sqrt{-(l_1 + e_3)}
\]
are all purely imaginary, and the same is true of the quantities
\begin{align*}
h_0 &= \frac{\sigma_1(w)}{\sigma(w)} + \frac{\sigma_2(w)}{\sigma(w)} + \frac{\sigma_3(w)}{\sigma(w)}, \\
h_1 &= (e_2 - e_3)\frac{\sigma_1(w)}{\sigma(w)} + (e_3 - e_1)\frac{\sigma_2(w)}{\sigma(w)} + (e_1 - e_2)\frac{\sigma_3(w)}{\sigma(w)}, \\
h_2 &= (e_2^2 - e_3^2)\frac{\sigma_1(w)}{\sigma(w)} + (e_3^2 - e_1^2)\frac{\sigma_2(w)}{\sigma(w)} + (e_1^2 - e_2^2)\frac{\sigma_3(w)}{\sigma(w)}.
\end{align*}
With the coefficient of $x^4$ in $R(x)$ being equal to $-1$, we shall therefore set in this case, denoting by $E$ the product $(e_2 - e_3)(e_3 - e_1)(e_1 - e_2)$,
\begin{align*}
p + qi &= x_1 = \varphi(u_1) = i\left(h_0 + \frac{2E}{h_1\wp u_1 + h_2}\right), \\
p - qi &= x_2 = \varphi(u_2) = i\left(h_0 + \frac{2E}{h_1\wp u_2 + h_2}\right),
\end{align*}
and one sees that $x_1$ and $x_2$ will be conjugate imaginary quantities, if $u_1$ and $u_2$ are. Let us set
\[
u_1 = \frac{u + vi}{2}, \quad u_2 = \frac{u - vi}{2},
\]
with $u$ and $v$ being real.
One then has
\[
s_1 = \wp(u_1 + u_2) = \wp(u), \quad s_2 = \wp(u_1 - u_2) = \wp(vi).
\]
It follows from these formulas that $s_1$ and $s_2$ are both real quantities, contained between the following limits
\[
(\infty \ldots s_1 \ldots e_1), \quad (e_3 \ldots s_2 \ldots -\infty).
\]

To express the quantities $p$ and $q$ by means of $s_1$ and $s_2$, I shall make use of the following identities.
One finds on page 50 of the cited work of Mr. SCHWARZ the following
\begin{align*}
&\sigma_\lambda(w)\sigma(u + v + w)\sigma(u - v) \\
&= \sigma(u + w)\sigma(u)\sigma_\lambda(v + w)\sigma_\lambda(v) - \sigma_\lambda(u + w)\sigma_\lambda(u)\sigma(v + w)\sigma(v), \\
&\sigma_\lambda(w)\sigma_\lambda(u + v + w)\sigma_\lambda(u - v) \\
&= \sigma_\lambda(u + w)\sigma_\lambda(u)\sigma_\lambda(v + w)\sigma_\lambda(v) - (e_\lambda - e_\mu)(e_\lambda - e_\nu)\sigma(u + w)\sigma(u)\sigma(v + w)\sigma(v).
\end{align*}

Setting $w = 0$ in these identities and giving $\lambda$ successively the values 1, 2, 3, one finds
\begin{align*}
\sigma(u + v)\sigma(u - v) &= \sigma^2u\sigma_1^2v - \sigma_1^2u\sigma^2v, \\
\sigma_1(u + v)\sigma_1(u - v) &= \sigma_1^2u\sigma_1^2v - (e_1 - e_2)(e_1 - e_3)\sigma^2u\sigma^2v, \\
\sigma_2(u + v)\sigma_2(u - v) &= \sigma_2^2u\sigma_2^2v + (e_1 - e_2)(e_2 - e_3)\sigma^2u\sigma^2v, \\
\sigma_3(u + v)\sigma_3(u - v) &= \sigma_3^2u\sigma_3^2v - (e_1 - e_3)(e_2 - e_3)\sigma^2u\sigma^2v,
\end{align*}
or, recognising that
\begin{align*}
\sigma_2^2u &= \sigma_1^2u + (e_1 - e_2)\sigma^2u, \\
\sigma_3^2u &= \sigma_1^2u + (e_1 - e_3)\sigma^2u,
\end{align*}
\begin{align*}
&\sigma_2(u + v)\sigma_2(u - v) = \sigma_1^2u\sigma_1^2v + (e_1 - e_2)(e_1 - e_3)\sigma^2u\sigma^2v + (e_1 - e_2)(\sigma^2u\sigma_1^2v + \sigma_1^2u\sigma^2v), \\
&\sigma_3(u + v)\sigma_3(u - v) = \sigma_1^2u\sigma_1^2v + (e_1 - e_2)(e_1 - e_3)\sigma^2u\sigma^2v + (e_1 - e_3)(\sigma^2u\sigma_1^2v + \sigma_1^2u\sigma^2v).
\end{align*}

It follows from these formulas that
\begin{align*}
&2(e_2 - e_3)(e_3 - e_1)(e_1 - e_2)\sigma^2u\sigma^2v = \\
&\qquad (e_2 - e_3)\sigma_1(u + v)\sigma_1(u - v) + (e_3 - e_1)\sigma_2(u + v)\sigma_2(u - v) + (e_1 - e_2)\sigma_3(u + v)\sigma_3(u - v), \\
&2(e_2 - e_3)\sigma_1^2u\sigma_1^2v = \\
&\qquad(e_2 - e_3)\sigma_1(u + v)\sigma_1(u - v) - (e_3 - e_1)\sigma_2(u + v)\sigma_2(u - v) - (e_1 - e_2)\sigma_3(u + v)\sigma_3(u - v), \\
&2(e_2 - e_3)\sigma^2u\sigma_1^2v = (e_2 - e_3)\sigma(u + v)\sigma(u - v) - \sigma_2(u + v)\sigma_2(u - v) + \sigma_3(u + v)\sigma_3(u - v), \\
&2(e_2 - e_3)\sigma_1^2u\sigma^2v = -(e_2 - e_3)\sigma(u + v)\sigma(u - v) - \sigma_2(u + v)\sigma_2(u - v) + \sigma_3(u + v)\sigma_3(u - v),
\end{align*}
from which one obtains by division
\begin{align*}
\frac{\sigma_1^2u}{\sigma^2u} &= (e_1 - e_2)(e_1 - e_3) \\
&\quad \times \frac{(e_2 - e_3)\sigma(u + v)\sigma(u - v) + \sigma_2(u + v)\sigma_2(u - v) - \sigma_3(u + v)\sigma_3(u - v)}{(e_2 - e_3)\sigma_1(u + v)\sigma_1(u - v) + (e_3 - e_1)\sigma_2(u + v)\sigma_2(u - v) + (e_1 - e_2)\sigma_3(u + v)\sigma_3(u - v)}
\end{align*}
and consequently, since $\wp u = \frac{\sigma_1^2u}{\sigma^2u} + e_1$,
\begin{align*}
\wp u = & -\frac{E\sigma(u + v)\sigma(u - v)}{(e_2 - e_3)\sigma_1(u + v)\sigma_1(u - v) + (e_3 - e_1)\sigma_2(u + v)\sigma_2(u - v) + (e_1 - e_2)\sigma_3(u + v)\sigma_3(u - v)} - \\
&\frac{(e_2^2 - e_3^2)\sigma_1(u + v)\sigma_1(u - v) + (e_3^2 - e_1^2)\sigma_2(u + v)\sigma_2(u - v) + (e_1^2 - e_2^2)\sigma_3(u + v)\sigma_3(u - v)}{(e_2 - e_3)\sigma_1(u + v)\sigma_1(u - v) + (e_3 - e_1)\sigma_2(u + v)\sigma_2(u - v) + (e_1 - e_2)\sigma_3(u + v)\sigma_3(u - v)}, \\
\wp v = & \frac{E\sigma(u + v)\sigma(u - v)}{(e_2 - e_3)\sigma_1(u + v)\sigma_1(u - v) + (e_3 - e_1)\sigma_2(u + v)\sigma_2(u - v) + (e_1 - e_2)\sigma_3(u + v)\sigma_3(u - v)} - \\
& \frac{(e_2^2 - e_3^2)\sigma_1(u + v)\sigma_1(u - v) + (e_3^2 - e_1^2)\sigma_2(u + v)\sigma_2(u - v) + (e_1^2 - e_2^2)\sigma_3(u + v)\sigma_3(u - v)}{(e_2 - e_3)\sigma_1(u + v)\sigma_1(u - v) + (e_3 - e_1)\sigma_2(u + v)\sigma_2(u - v) + (e_1 - e_2)\sigma_3(u + v)\sigma_3(u - v)}.
\end{align*}
Setting therefore
\[
P_a = \frac{\sigma_a(u_1 + u_2)\sigma_a(u_1 - u_2)}{\sigma(u_1 + u_2)\sigma(u_1 - u_2)} = \sqrt{\left(s_1-e_a\right)\left(s_2-e_a\right)}, \quad (a = 1, 2, 3)
\]
one has
\begin{align*}
\wp(u_1) + \wp(u_2) &= -2\frac{(e_2^2 - e_3^2)P_1 + (e_3^2 - e_1^2)P_2 + (e_1^2 - e_2^2)P_3}{(e_2 - e_3)P_1 + (e_3 - e_1)P_2 + (e_1 - e_2)P_3}, \\
\wp(u_1) - \wp(u_2) &= \frac{-2E}{(e_2 - e_3)P_1 + (e_3 - e_1)P_2 + (e_1 - e_2)P_3}, \\
\wp(u_1)\wp(u_2) &= -\frac{(e_2 - e_3)(e_1^2 + e_2e_3)P_1 + (e_3 - e_1)(e_2^2 + e_3e_1)P_2 + (e_1 - e_2)(e_3^2 + e_1e_2)P_3}{(e_2 - e_3)P_1 + (e_3 - e_1)P_2 + (e_1 - e_2)P_3}.
\end{align*}

Substituting these expressions into the formulas
\begin{align*}
p &= \frac{x_1 + x_2}{2} = i\left(h_0 + E\frac{h_1(\wp u_1 + \wp u_2) + 2h_2}{h_1^2\wp u_1\wp u_2 + h_1h_2(\wp u_1 + \wp u_2) + h_2^2}\right), \\
q &= \frac{x_1 - x_2}{2i} = E\frac{h_1(\wp u_1 - \wp u_2)}{h_1^2\wp u_1\wp u_2 + h_1h_2(\wp u_1 + \wp u_2) + h_2^2},
\end{align*}
one finds, after some calculations, that
\begin{align*}
p &= -i\frac{L_1P_1 + M_1P_2 + N_1P_3}{LP_1 + MP_2 + NP_3}, \\
q &= \frac{E}{LP_1 + MP_2 + NP_3},
\end{align*}
where I have set
\begin{align*}
L &= (e_2 - e_3)\frac{\sigma_1(w)}{\sigma(w)} = i(e_2 - e_3)\sqrt{l_1 + e_1}, \\
M &= (e_3 - e_1)\frac{\sigma_2(w)}{\sigma(w)} = i(e_3 - e_1)\sqrt{l_1 + e_2}, \\
N &= (e_1 - e_2)\frac{\sigma_3(w)}{\sigma(w)} = i(e_1 - e_2)\sqrt{l_1 + e_3}, \\
L_1 &= (e_2 - e_3)\frac{\sigma_2(w)\sigma_3(w)}{\sigma^2(w)} = (e_2 - e_3)\sqrt{(l_1 + e_2)(l_1 + e_3)}, \\
M_1 &= (e_3 - e_1)\frac{\sigma_3(w)\sigma_1(w)}{\sigma^2(w)} = (e_3 - e_1)\sqrt{(l_1 + e_3)(l_1 + e_1)}, \\
N_1 &= (e_1 - e_2)\frac{\sigma_1(w)\sigma_2(w)}{\sigma^2(w)} = (e_1 - e_2)\sqrt{(l_1 + e_1)(l_1 + e_2)}.
\end{align*}

The quantities $l_1 + e_1$, $l_1 + e_2$, $l_1 + e_3$ are all positive.
The quantities $s_1 - e_1$, $s_1 - e_2$, $s_1 - e_3$ are also positive.
The quantities $s_2 - e_1$, $s_2 - e_2$, $s_2 - e_3$ are on the contrary all negative.
$L$, $M$, $N$, $P_1$, $P_2$, $P_3$ are therefore imaginary.
$L_1$, $M_1$, $N_1$ are real.
One therefore sees that $p$ and $q$ have real values.

To calculate the value of $r$ I make use of the equation
\[
2\frac{dp}{dt} = qr.
\]

In general, if one sets
\[
R(s) = R_0(s - a_0)(s - a_1)(s - a_2)(s - a_3)(s - a_4),
\]
\begin{align*}
du_1 &= \frac{s_1ds_1}{\sqrt{R(s_1)}} + \frac{s_2ds_2}{\sqrt{R(s_2)}}, \quad &ds_1 &= \frac{\sqrt{R(s_1)}}{s_1 - s_2}(du_1 - s_2du_2), \\
du_2 &= \frac{ds_1}{\sqrt{R(s_1)}} + \frac{ds_2}{\sqrt{R(s_2)}}, \quad &ds_2 &= \frac{\sqrt{R(s_2)}}{s_2 - s_1}(du_1 - s_1du_2),
\end{align*}
\begin{align*}
P_a &= \sqrt{c_a(s_1 - a_a)(s_2 - a_a)} \\
P_{a\beta} &= \frac{c_{a\beta}P_aP_\beta}{s_1 - s_2}\left\{\frac{\sqrt{R(s_1)}}{(s_1 - a_a)(s_1 - a_\beta)} - \frac{\sqrt{R(s_2)}}{(s_2 - a_a)(s_2 - a_\beta)}\right\}, \quad \begin{pmatrix} a = 0, \ldots, 4 \\ \beta = 0, \ldots, 4 \end{pmatrix}
\end{align*}
with $c_\alpha$ and $c_{\alpha\beta}$ denoting arbitrary constants, one finds, after some calculations, that
\begin{align*}
\frac{\partial P_a}{\partial u_1} &= \frac{1}{2(a_\beta - a_\gamma)}\left\{\frac{P_\gamma P_{a\gamma}}{c_\gamma c_{a\gamma}} - \frac{P_\beta P_{a\beta}}{c_\beta c_{a\beta}}\right\}, \\
\frac{\partial P_a}{\partial u_2} &= \frac{1}{2(a_\beta - a_\gamma)}\left\{\frac{a_\gamma P_\beta P_{a\beta}}{c_\beta c_{a\beta}} - \frac{a_\beta P_\gamma P_{a\gamma}}{c_\gamma c_{a\gamma}}\right\}, \\
\frac{\partial P_{a\beta}}{\partial u_1} &= \frac{1}{2}R_0 c_{a\beta} P_a P_\beta, \\
\frac{\partial P_{a\beta}}{\partial u_2} &= -\frac{1}{2}R_0 a_\gamma c_{a\beta} P_a P_\beta - \frac{1}{2}\frac{c_{a\beta}P_{a\gamma}P_{\beta\gamma}}{c_\gamma c_{a\gamma}c_{\beta\gamma}}.
\end{align*}

In our case one has
\begin{align*}
R(s) &= -4(s - e_1)(s - e_2)(s - e_3)(s - k_1)(s - k_2), \quad R_0 = -4, \\
dt &= du_1, \quad 0 = du_2.
\end{align*}

One therefore finds that
\[
r = \frac{2}{q}\frac{dp}{dt} = -i\frac{LP_{23} + MP_{31} + NP_{12}}{LP_1 + MP_2 + NP_3},
\]
\begin{align*}
P_{23} &= \frac{i}{s_1 - s_2}\left[\sqrt{(s_1 - e_2)(s_1 - e_3)(s_1 - k_1)(s_1 - k_2)} - \sqrt{(s_2 - e_2)(s_2 - e_3)(s_2 - k_1)(s_2 - k_2)}\right], \\
P_{13} &= i\frac{\sqrt{(s_1 - e_2)(s_1 - k_1)(s_1 - k_2)(s_2 - e_1)(s_2 - e_3)}}{s_1\left(s_1 - s_2\right)} - i\frac{\sqrt{(s_1 - e_3)(s_1 - e_1)(s_2 - e_2)(s_2 - k_1)(s_2 - k_2)}}{s_2\left(s_1 - s_2\right)}, \\
P_{12} &= i\frac{\sqrt{(s_1 - e_3)(s_1 - k_1)(s_1 - k_2)(s_2 - e_1)(s_2 - e_2)}}{s_1 - s_2} - i\frac{\sqrt{(s_2 - k_2)(s_1 - e_1)(s_1 - e_2)(s_2 - e_3)(s_2 - k_1)}}{s_1 - s_2}.
\end{align*}

As we have seen above, for real values of $t$, $s_1$ and $s_2$ must be contained between the following limits
\[
(\infty \ldots s_1 \ldots e_1), \quad (e_3 \ldots s_2 \ldots -\infty).
\]
The quantities
\begin{gather*}
(s_1 - e_1)(s_2 - e_2)(s_2 - e_3) \\
(s_1 - e_2)(s_2 - e_3)(s_2 - e_1) \\
(s_1 - e_3)(s_2 - e_1)(s_2 - e_2)
\end{gather*}
are therefore positive, and the following are negative
\begin{gather*}
(s_1 - e_2)(s_1 - e_3)(s_2 - e_1) \\
(s_1 - e_3)(s_1 - e_1)(s_2 - e_2) \\
(s_1 - e_1)(s_1 - e_2)(s_2 - e_3).
\end{gather*}

For $r$ to be real for real values of $t$, it is necessary that $P_{23}$, $P_{13}$, $P_{12}$ be real; for this to take place, it is necessary that the product
\[
(s_1 - k_1)(s_1 - k_2)
\]
is negative, and that the following product is positive
\[
(s_2 - k_1)(s_2 - k_2).
\]
The constants of integration must therefore be chosen in such a manner to satisfy the inequality
\[
k_1 > e_1 > k_2.
\]
therefore $s_1$ and $s_2$ must then be contained between the limits
\[
(k_1 \ldots s_1 \ldots e_1), \quad (k_2 \ldots s_2 \ldots -\infty).
\]

It remains for us to calculate the values of $\gamma$, $\gamma'$, $\gamma''$.
One obtains $\gamma''$ from the equation
\[
2\frac{dq}{dt} = -rp - c_0\gamma''
\]
or
\[
c_0\gamma'' = -\left(2\frac{dq}{dt} + rp\right),
\]
and one finds with the aid of the aforementioned differentiation formulas
\[
c_0\gamma'' = \frac{L_1P_{23} + M_1P_{13} + N_1P_{12}}{LP_{23} + MP_{13} + NP_{12}},
\]
the coefficients $L$, $M$, $N$, $L_1$, $M_1$, $N_1$ having the same meanings as before.

The value of $\gamma'$ can be calculated with the aid of the equation
\[
\frac{dr}{dt} = c_0\gamma'.
\]

One finds
\begin{align*}
c_0\gamma' &= -i\frac{(LP_1 + MP_2 + NP_3)\left(L\frac{d}{dt}P_{23} + M\frac{d}{dt}P_{13} + N\frac{d}{dt}P_{12}\right)}{(LP_1 + MP_2 + NP_3)^2} \\
&\quad + i\frac{(LP_{23} + MP_{13} + NP_{12})\left(L\frac{d}{dt}P_1 + M\frac{d}{dt}P_2 + N\frac{d}{dt}P_3\right)}{(LP_1 + MP_2 + NP_3)^2}.
\end{align*}

According to the aforementioned differentiation formulas, one finds
\begin{align*}
&(LP_1 + MP_2 + NP_3)\left(L\frac{d}{dt}P_{23} + M\frac{d}{dt}P_{13} + N\frac{d}{dt}P_{12}\right) \\
&\quad - (LP_{23} + MP_{13} + NP_{12})\left(L\frac{d}{dt}P_1 + M\frac{d}{dt}P_2 + N\frac{d}{dt}P_3\right) \\
&= R_0(LP_1 + MP_2 + NP_3)(LP_2P_3 + MP_3P_1 + NP_1P_2) \\
&\quad - (LP_{23} + MP_{31} + NP_{12})\left(L\frac{P_3P_{13} - P_2P_{12}}{e_2 - e_3} + M\frac{P_1P_{12} - P_3P_{23}}{e_3 - e_1} + N\frac{P_2P_{23} - P_1P_{13}}{e_1 - e_2}\right) \\
&= R_0(L^2 + M^2 + N^2)P_1P_2P_3 - \frac{L^2}{e_2 - e_3}P_{23}(P_3P_{13} - P_2P_{12}) \\
&\quad - \frac{M^2}{e_3 - e_1}P_{13}(P_1P_{12} - P_3P_{23}) - \frac{N^2}{e_1 - e_2}P_{12}(P_2P_{23} - P_1P_{13}) \\
&\quad + MN\left(R_0(P_2^2 + P_3^2)P_1 - P_{13}\frac{P_2P_{23} - P_1P_{13}}{e_1 - e_2} - P_{12}\frac{P_1P_{12} - P_3P_{23}}{e_3 - e_1}\right) \\
&\quad + NL\left(R_0(P_3^2 + P_1^2)P_2 - P_{12}\frac{P_3P_{13} - P_2P_{12}}{e_2 - e_3} - P_{23}\frac{P_2P_{23} - P_1P_{13}}{e_1 - e_2}\right) \\
&\quad + LM\left(R_0(P_1^2 + P_2^2)P_3 - P_{23}\frac{P_1P_{12} - P_3P_{23}}{e_3 - e_1} - P_{13}\frac{P_3P_{13} - P_2P_{12}}{e_2 - e_3}\right).
\end{align*}

This expression for the numerator of $\gamma'$ can be somewhat simplified in the following manner.
Denoting by $\alpha$, $\beta$, $\gamma$, $\delta$, $\epsilon$ the numbers 1, 2, 3, 4, 5 arranged in any order and setting, as we have done,
\[
P_\alpha = \sqrt{(s_1 - a_\alpha)(s_2 - a_\alpha)},
\]
\begin{align*}
P_{\alpha\beta} &= \frac{\sqrt{R_0(s_1 - a_\gamma)(s_1 - a_\delta)(s_1 - a_\epsilon)(s_2 - a_a)(s_2 - a_\beta)}}{s_1 - s_2} \\
&\quad - \frac{\sqrt{R_0(s_1 - a_a)(s_1 - a_\beta)(s_2 - a_\gamma)(s_2 - a_\delta)(s_2 - a_\epsilon)}}{s_1 - s_2},
\end{align*}
one easily finds the following relations:
\begin{align*}
&R_0P_\alpha P_\beta P_\gamma - P_{\beta\gamma}\frac{P_\gamma P_{\alpha\gamma} - P_\beta P_{\alpha\beta}}{a_\beta - a_\gamma} = P_{\alpha\delta}P_{\alpha\epsilon}, \\
&R_0P_\alpha P_\beta P_\gamma - P_{\alpha\gamma}\frac{P_\alpha P_{\alpha\beta} - P_\gamma P_{\beta\gamma}}{a_\gamma - a_\alpha} = P_{\beta\delta}P_{\beta\epsilon}, \\
&R_0P_\alpha P_\beta P_\gamma - P_{\alpha\beta}\frac{P_\beta P_{\beta\gamma} - P_\alpha P_{\alpha\gamma}}{a_\alpha - a_\beta} = P_{\gamma\delta}P_{\gamma\epsilon}, \\
&R_0(P_\beta^2 + P_\gamma^2)P_\alpha - P_{\alpha\gamma}\frac{P_\beta P_{\beta\gamma} - P_\alpha P_{\alpha\gamma}}{a_\alpha - a_\beta} - P_{\alpha\beta}\frac{P_\alpha P_{\alpha\beta} - P_\gamma P_{\beta\gamma}}{a_\gamma - a_\alpha} \\
&\qquad = P_{\beta\delta}P_{\gamma\epsilon} + P_{\beta\epsilon}P_{\gamma\delta} + R_0(a_\beta - a_\gamma)^2P_\alpha, \\
&R_0(P_\gamma^2 + P_\alpha^2)P_\beta - P_{\alpha\beta}\frac{P_\gamma P_{\alpha\gamma} - P_\beta P_{\alpha\beta}}{a_\beta - a_\gamma} - P_{\beta\gamma}\frac{P_\beta P_{\beta\gamma} - P_\alpha P_{\alpha\gamma}}{a_\alpha - a_\beta} \\
&\qquad = P_{\gamma\delta}P_{\alpha\epsilon} + P_{\alpha\delta}P_{\gamma\epsilon} + R_0(a_\gamma - a_\alpha)^2P_\beta, \\
&R_0(P_\alpha^2 + P_\beta^2)P_\gamma - P_{\beta\gamma}\frac{P_\alpha P_{\alpha\beta} - P_\gamma P_{\beta\gamma}}{a_\gamma - a_\alpha} - P_{\alpha\gamma}\frac{P_\gamma P_{\alpha\gamma} - P_\beta P_{\alpha\beta}}{a_\beta - a_\gamma} \\
&\qquad = P_{\alpha\delta}P_{\beta\epsilon} + P_{\alpha\epsilon}P_{\beta\delta} + R_0(a_\alpha - a_\beta)^2P_\gamma.
\end{align*}
If one therefore sets
\begin{align*}
&a_1 = e_1, \\
&a_2 = e_2, \\
&a_3 = e_3, \\
&k_1 = a_4, \\
&k_2 = a_5,
\end{align*}
one finds for $c_0\gamma'$ the following expression
\begin{align*}
c_0\gamma' &= \frac{i}{2(LP_1 + MP_2 + NP_3)^2}\left[L^2P_{14}P_{15} + M^2P_{24}P_{25} + N^2P_{34}P_{35}\right. \\
&\quad + MN[P_{24}P_{35} + P_{25}P_{34} + (e_2 - e_3)^2P_1] \\
&\quad + NL[P_{34}P_{15} + P_{35}P_{14} + (e_3 - e_1)^2P_2] \\
&\quad \left. + LM[P_{14}P_{25} + P_{15}P_{24} + (e_1 - e_2)^2P_3]\right].
\end{align*}

Finally one finds the value of $\gamma$, expressed with the aid of $s_1$ and $s_2$, with the aid of the equation
\[
2(p^2 + q^2) + r^2 = 2c_0\gamma + 6l_1.
\]
Setting
\begin{align*}
L_2 &= (e_2^2 - e_3^2)\frac{\sigma_1(w)}{\sigma(w)} = i(e_2^2 - e_3^2)\sqrt{l_1 + e_1}, \\
M_2 &= (e_3^2 - e_1^2)\frac{\sigma_2(w)}{\sigma(w)} = i(e_3^2 - e_1^2)\sqrt{l_1 + e_2}, \\
N_2 &= (e_1^2 - e_2^2)\frac{\sigma_3(w)}{\sigma(w)} = i(e_1^2 - e_2^2)\sqrt{l_1 + e_3},
\end{align*}
one finds
\[
2c_0\gamma = -4l_1 + \frac{L_2P_1 + M_2P_2 + N_2P_3}{LP_1 + MP_2 + NP_3} - \left(\frac{LP_23 + MP_13 + NP_12}{LP_1 + MP_2 + NP_3}\right)^2.
\]

\section*{§ 6.}
\setcounter{equation}{0}

When the six quantities $p, q, r, \gamma, \gamma', \gamma''$ are expressed as functions of $s_1$ and of $s_2$, one obtains without difficulty the expression of each rational symmetric function of these two last quantities, as a function of time, with the aid of general formulas.

I borrow the following definitions, introduced in the analysis by Mr. WEIERSTRASS, from the memoir of Mr. KÖNIGSBERGER \textit{Zur Transformation der Abelschen Functionen} (Journal für die reine und angewandte Mathematik, Bd. 64).

Let
\[
R(x) = A_0(x - a_0)(x - a_1)\ldots(x - a_{2\rho})
\]
be an entire function of $x$ of degree $2\rho + 1$; suppose that $A_0$, $a_0$, $a_1, \ldots, a_{2\rho}$ are all real quantities and, if $A_0 > 0$,
\[
a_0 > a_1 > \ldots > a_{2\rho},
\]
but if $A_0 < 0$
\[
a_0 < a_1 < \ldots < a_{2\rho}.
\]
Let $u_1, \ldots, u_\rho$ be $\rho$ variables, linked to the $\rho$ variables $x_1, \ldots, x_\rho$ by the $\rho$ following equations, in which $F_1(x), \ldots, F_\rho(x)$ denote entire functions of $x$ of degree less than $\rho$,
\begin{align*}
u_1 &= \int_{a_1}^{x_1} \frac{F_1(x)dx}{\sqrt{R(x)}} + \ldots + \int_{a_{2\rho-1}}^{x_\rho} \frac{F_1(x)dx}{\sqrt{R(x)}}, \\
&\vdots \\
u_\rho &= \int_{a_1}^{x_1} \frac{F_\rho(x)dx}{\sqrt{R(x)}} + \ldots + \int_{a_{2\rho-1}}^{x_\rho} \frac{F_\rho(x)dx}{\sqrt{R(x)}}.
\end{align*}
Let us set
\begin{align*}
K_{\alpha\beta} &= \int_{a_{2\beta-1}}^{a_{2\beta}} \frac{F_\alpha(x)dx}{\sqrt{R(x)}}, \\
i\bar{K}_{\alpha\beta} &= \int_{a_{2\beta-2}}^{a_{2\beta-1}} \frac{F_\alpha(x)dx}{\sqrt{R(x)}},
\end{align*}
\[
iK'_{\mu\nu} = i\bar{K}_{\mu 1} + i\bar{K}_{\mu 2} + \ldots + i\bar{K}_{\mu\nu},
\]
agreeing to define the square root in each of these integrals by the formula
\[
\sqrt{R(x)} = A_0^{\rho+1}\left(\frac{x - a_0}{A_0}\right)^{\frac{1}{2}}\left(\frac{x - a_1}{A_0}\right)^{\frac{1}{2}}\ldots\left(\frac{x - a_{2\rho}}{A_0}\right)^{\frac{1}{2}}.
\]
Let us define $\rho$ new variables $v_1, \ldots, v_\rho$ by the equations
\begin{align*}
u_1 &= 2K_{11}v_1 + \ldots + 2K_{1\rho}v_\rho, \\
&\vdots \\
u_\rho &= 2K_{\rho 1}v_1 + \ldots + 2K_{\rho\rho}v_\rho.
\end{align*}
Suppose that these equations, solved with respect to $v_1, \ldots, v_\rho$, give us
\begin{align*}
v_1 &= G_{11}u_1 + \ldots + G_{\rho 1}u_\rho, \\
&\vdots \\
v_\rho &= G_{1\rho}u_1 + \ldots + G_{\rho\rho}u_\rho.
\end{align*}

Let us set
\[
\tau_{\alpha\beta} = 2i(G_{1\alpha}K'_{1\beta} + \ldots + G_{\rho\alpha}K'_{\rho\beta}),
\]
from which it follows that $\tau_{\alpha\beta} = \tau_{\beta\alpha}$, and let us define the function $\vartheta(v_1\ldots v_\rho)$ as the sum of the infinite series
\[
\vartheta(v_1\ldots v_\rho) = \sum e^{\left\{\nu_1(2v_1 + \nu_1\tau_{11} + \ldots + \nu_\rho\tau_{1\rho}) + \ldots + \nu_\rho(2v_\rho + \nu_1\tau_{\rho 1} + \ldots + \nu_\rho\tau_{\rho\rho})\right\}\pi i};
\]
where the summation, indicated by the sign $\sum$, must be performed in such a manner that each of the $\rho$ quantities $\nu_1, \ldots, \nu_\rho$ traverses, independently of the others, the entire series of integers from $-\infty$ to $+\infty$.

This function $\vartheta(v_1, \ldots, v_\rho)$ satisfies the two following equations:
\[
\vartheta(v_1 + p_1, \ldots, v_\rho + p_\rho) = \vartheta(v_1\ldots v_\rho),
\]
with $p_1\ldots p_\rho$ denoting any integers, and
\[
\vartheta(v_1 + \tau_{1\alpha}, \ldots, v_\rho + \tau_{\rho\alpha}) = \vartheta(v_1\ldots v_\rho)e^{-\pi i(2v_\alpha + \tau_{\alpha\alpha})}.
\]
Conversely, each continuous function of $v_1\ldots v_\rho$ satisfying these two equations for all systems of values $v_1\ldots v_\rho$, must necessarily be equal to $\vartheta(v_1\ldots v_\rho)$, multiplied by a constant.

Let us denote by $n_1\ldots n_\rho$ arbitrary constants, let us set
\[
\tau_\alpha = n_1\tau_{\alpha 1} + n_2\tau_{\alpha 2} + \ldots + n_\rho\tau_{\alpha \rho}
\]
and let us define a new function $\vartheta(v_1\ldots v_\rho \mid n_1\ldots n_\rho)$ by the equation
\[
\vartheta(v_1\ldots v_\rho \mid n_1\ldots n_\rho) = \vartheta(v_1 + \tau_1\ldots v_\rho + \tau_\rho)e^{\pi i\sum_\alpha n_\alpha(2v_\alpha + \tau_\alpha)}.
\]

The quantities $n_1\ldots n_\rho$ are called the parameters of this new function $\vartheta$. If $n'_1\ldots n'_\rho$ denote another system of $\rho$ constants, and one sets
\[
\tau'_\alpha = n'_1\tau'_{1\alpha} + \ldots + n'_\rho\tau'_{\rho\alpha},
\]
one has
\[
\vartheta(v_1 + \tau'_{1}, \ldots, v_\rho + \tau'_{\rho}\mid n_1\ldots n_\rho) = \vartheta(v_1\ldots v_\rho \mid n_1 + n'_1\ldots n_\rho + n'_\rho)e^{-\pi i\sum_\alpha n'_\alpha(2v_\alpha + \tau'_\alpha)}.
\]

If the quantities $n'_1\ldots n'_\rho$, $m_1\ldots m_\rho$ are integers, one has
\begin{align*}
\vartheta(v_1\ldots v_\rho \mid n_1 + n'_1\ldots n_\rho + n'_\rho) &= \vartheta(v_1\ldots v_\rho \mid n_1\ldots n_\rho), \\
\vartheta(v_1 + m_1\ldots v_\rho + m_\rho \mid n_1\ldots n_\rho) &= e^{2\pi i\sum_{\alpha} m_\alpha n_\alpha}\vartheta(v_1\ldots v_\rho \mid n_1\ldots n_\rho), \\
\vartheta(v_1 + \tau'_1\ldots v_\rho + \tau'_\rho \mid n_1\ldots n_\rho) &= e^{-\pi i\sum n_\alpha(2v_\alpha + \tau'_\alpha)}\vartheta(v_1\ldots v_\rho \mid n_1\ldots n_\rho).
\end{align*}

Let us set moreover
\[
\vartheta(v_1\ldots v_\rho)_\lambda = \vartheta\left(v_1 + \frac{1}{2}m^{\lambda}_{1}\ldots v_\rho + \frac{1}{2}m^{\lambda}_{\rho} \mid \frac{1}{2}n^{\lambda}_{1}\ldots\frac{1}{2}n^{\lambda}_{\rho}\right) \quad (\lambda = 0, 1, 2, \ldots, 2\rho)
\]
where the integers $m^{\lambda}_{1}\ldots m^{\lambda}_{\rho}$, $n^{\lambda}_{1}\ldots n^{\lambda}_{\rho}$ are defined by the equation
\[
\int_{\infty}^{a_\lambda} \frac{F_\alpha(x)dx}{\sqrt{R(x)}} = m^{\lambda}_{1}K_{\alpha 1} + \ldots + m^{\lambda}_{\rho}K_{\alpha\rho} + i(n^{\lambda}_{1}K'_{\alpha 1} + n^{\lambda}_{2}K'_{\alpha 2} + \ldots + n^{\lambda}_{\rho}K'_{\alpha\rho}).
\]

(One easily verifies that each of the integers $m^{\lambda}_{1}\ldots m^{\lambda}_{\rho}$ is equal to 0 or $-1$; and each of the integers $n^{\lambda}_{1}\ldots n^{\lambda}_{\rho}$ to 0 or $+1$).

Denoting next by $\lambda$ and $\mu$ two different arbitrary numbers from the series $0, 1, \ldots, 2\rho$, let us define the integers $m_1^\nu\ldots m_\rho^\nu$, $n_1^\nu\ldots n_\rho^\nu$ by the congruences
\[
\begin{cases}
m_a^\nu \equiv m_a^\lambda + m_a^\mu \\
n_a^\nu \equiv n_a^\lambda + n_a^\mu
\end{cases} \pmod{2}
\]
and moreover by the condition that each of the numbers $m_\lambda^\nu$ must be equal to 0 or to $-1$, and each of the numbers $n_\lambda^\nu$ to 0 or to $+1$, and let us set
\[
\vartheta(v_1\ldots v_\rho)_{\lambda\mu} = \vartheta\left(v_1 + \frac{1}{2}m_1^\nu\ldots v_\rho + \frac{1}{2}m_\rho^\nu \,\middle|\, \frac{1}{2}n_1^\nu\ldots \frac{1}{2}n_\rho^\nu\right).
\]

The relation between the variables $v_1\ldots v_\rho$ and $x_1\ldots x_\rho$ can then be expressed in the following manner.

Let us set
\[
\varphi(x) = (x - x_1)\ldots(x - x_\rho)
\]
and denote by $\epsilon = \pm 1$, according to $A_0 \gtrless 0$; one then has
\begin{align*}
\frac{\sqrt{\epsilon^{\rho}(-1)^{\alpha}\varphi(a_{2\alpha})}}{\sqrt[4]{R'(a_{2\alpha})}} &= \frac{\vartheta(v_1\ldots v_\rho)_{2\alpha}}{\vartheta(v_1\ldots v_\rho)}, \\
\frac{\sqrt{\epsilon^{\rho}(-1)^{\alpha-1}\varphi(a_{2\alpha-1})}}{\sqrt[4]{R'(a_{2\alpha-1})}} &= \frac{\vartheta(v_1\ldots v_\rho)_{2\alpha-1}}{\vartheta(v_1\ldots v_\rho)}, \\
A_0\sqrt{\frac{\pm(a_\lambda - a_\mu)}{A_0}}\sum_{a=1}^{\rho}\left\{\frac{\sqrt{R(x_a)}}{(x_a - a_\lambda)(x_a - a_\mu)\varphi'(x_a)}\right\} &= \frac{\vartheta(v_1\ldots v_\rho)\vartheta(v_1\ldots v_\rho)_{\lambda\mu}}{\vartheta(v_1\ldots v_\rho)_\lambda\vartheta(v_1\ldots v_\rho)_\mu}.
\end{align*}

We have supposed that all the roots of the equation $R(x) = 0$ are real. Mr. HENOCH has shown in his inaugural dissertation (Berlin 1867) how these formulas can be generalized, for the case where the roots of the equation $R(x) = 0$ are imaginary.

\section*{§ 7.}
\setcounter{equation}{0}

To be able to apply the formulas of the preceding section to the case which occupies us, it is necessary to arrange the five real quantities $k_1$, $k_2$, $e_1$, $e_2$, $e_3$ in order of magnitude.

Several cases can present themselves here. Lack of time prevents me from examining all of them in detail. I shall content myself with carrying out all the calculations for the case where one has the following inequalities between the constants of integration
\[
l_1 > k > c_0 > 0, \quad l^2 < \frac{3l_1 - k}{2}.
\]

In this case, if one writes
\[
4s^3 - g_2s - g_3 = 4(s + l_1)\left(s - \frac{l_1 + \sqrt{k^2 - c_0^2}}{2}\right)\left(s - \frac{l_1 - \sqrt{k^2 - c_0^2}}{2}\right) - l_0^2,
\]
one easily verifies that the five quantities $e_1, e_2, e_3, k_1, k_2$ satisfy the inequalities
\[
\frac{l_1 + k}{2} > e_1 > e_2 > \frac{l_1 - k}{2} > e_3.
\]

If one denotes by $a_0, a_1, a_2, a_3, a_4$ the five roots of the equation $R(s) = 0$, arranged in order of magnitude
\[
a_0 < a_1 < a_2 < a_3 < a_4,
\]
it is therefore necessary to set
\[
\frac{l_1 + k}{2} = a_0, \quad \frac{l_1 - k}{2} = a_1, \quad e_1 = a_3, \quad e_2 = a_2, \quad e_3 = a_0.
\]

If one then sets
\[
dt = du_1 = \frac{s_1ds_1}{\sqrt{R(s_1)}} + \frac{s_2ds_2}{\sqrt{R(s_2)}},
\]
\[
0 = du_2 = \frac{ds_1}{\sqrt{R(s_1)}} + \frac{ds_2}{\sqrt{R(s_2)}},
\]
and one determines the quantities $K_{\alpha\beta}$, $K'_{\alpha\beta}$, $v_1, \ldots, v_\rho$ as it has been shown in the preceding section; $v_1$, $v_2$ will be linear functions of time, and one will have, denoting by $c_\lambda$ the value that $\vartheta(v_1, v_2)_\lambda$ takes for $v_1 = 0$, $v_2 = 0$, the following relations:
\begin{align*}
\frac{a_3 - a_2}{a_4 - a_1} &= \frac{e_1 - e_2}{k} = \frac{c_2^2c_{03}^2c_{34}^2}{c_4^2c_{01}^2c_{12}^2}, & \frac{l_1 + e_1}{k} &= \frac{a_4 + a_1 - a_2 - a_0}{a_4 - a_1} = \frac{c_5^2c_{23}^2}{c_4^2c_{01}^2} - \frac{c_2^2c_{14}^2}{c_4^2c_{12}^2}, \\
\frac{a_3 - a_0}{a_4 - a_1} &= \frac{e_1 - e_3}{k} = \frac{c_0^2c_{23}^2c_{34}^2}{c_4^2c_{01}^2c_{12}^2}, & \frac{l_1 + e_2}{k} &= \frac{a_4 + a_1 - a_0 - a_3}{a_4 - a_1} = \frac{c_5^2c_{23}^2}{c_4^2c_{01}^2} - \frac{c_0^2c_2^2}{c_{01}^2c_{12}^2}, \\
\frac{a_2 - a_0}{a_4 - a_1} &= \frac{e_2 - e_3}{k} = \frac{c_5^2c_{14}^2c_{34}^2}{c_4^2c_{01}^2c_{12}^2}, & \frac{l_1 + e_3}{k} &= \frac{a_4 + a_1 - a_3 - a_2}{a_4 - a_2} = \frac{c_5^2c_{03}^2}{c_4^2c_{12}^2} - \frac{c_0^2c_2^2}{c_{01}^2c_{12}^2},
\end{align*}
and, setting
\[
C = (a_3 - a_1)^{\frac{3}{2}}\frac{c_{01}c_{03}c_{12}c_{23}c_{14}c_{34}}{c_0^2c_2^2c_4^2},
\]
one has \hypertarget{TN3-ref}{}\tnref{TN3}:
\begin{align*}
P_1 &= \sqrt{(s_1 - e_1)(s_2 - e_1)} = \sqrt{(s_1 - a_3)(s_2 - a_3)} = -i\frac{C}{\sqrt{a_3 - a_1}}\frac{c_0c_2c_4}{c_{01}c_{12}c_{14}}\frac{\vartheta_3(v_1v_2)}{\vartheta_5(v_1v_2)}, \\
P_2 &= \sqrt{(s_1 - e_2)(s_2 - e_2)} = \sqrt{(s_1 - a_2)(s_2 - a_2)} = -i\frac{C}{\sqrt{a_3 - a_1}}\frac{c_5c_2}{c_{12}c_{23}}\frac{\vartheta_2(v_1v_2)}{\vartheta_5(v_1v_2)}, \\
P_3 &= \sqrt{(s_1 - e_3)(s_2 - e_3)} = \sqrt{(s_1 - a_0)(s_2 - a_0)} = \frac{C}{\sqrt{a_3 - a_1}}\frac{c_5c_0}{c_{01}c_{03}}\frac{\vartheta_0(v_1v_2)}{\vartheta_5(v_1v_2)}, \\
P_4 &= \sqrt{(s_1 - k_1)(s_2 - k_1)} = \sqrt{(s_1 - a_4)(s_2 - a_4)} = -\frac{C}{\sqrt{a_3 - a_1}}\frac{c_5c_4}{c_{14}c_{34}}\frac{\vartheta_4(v_1v_2)}{\vartheta_5(v_1v_2)}, \\
P_5 &= \sqrt{(s_1 - k_2)(s_2 - k_2)} = \sqrt{(s_1 - a_1)(s_2 - a_1)} = \frac{C}{\sqrt{a_3 - a_1}}\frac{c_0c_2c_4}{c_{03}c_{23}c_{34}}\frac{\vartheta_1(v_1v_2)}{\vartheta_5(v_1v_2)}.
\end{align*}

\begin{align*}
P_{12} &= \frac{\sqrt{(s_1 - e_3)(s_1 - k_1)(s_1 - k_2)(s_2 - e_1)(s_2 - e_2)} - \sqrt{(s_1 - e_1)(s_1 - e_2)(s_2 - e_3)(s_2 - k_1)(s_2 - k_2)}}{s_1 - s_2} \\
&= iC\frac{c_5}{c_{12}}\frac{\vartheta_{23}(v_1v_2)}{\vartheta_5(v_1v_2)},
\end{align*}
\begin{align*}
P_{13} &= \frac{\sqrt{(s_1 - e_2)(s_1 - k_1)(s_1 - k_2)(s_2 - e_1)(s_2 - e_3)} - \sqrt{(s_1 - e_1)(s_1 - e_3)(s_2 - e_2)(s_2 - k_1)(s_2 - k_2)}}{s_1 - s_2} \\
&= -C\frac{c_5}{c_{01}}\frac{\vartheta_{03}(v_1v_2)}{\vartheta_5(v_1v_2)},
\end{align*}
\begin{align*}
P_{14} &= \frac{\sqrt{(s_1 - e_2)(s_1 - e_3)(s_1 - k_2)(s_2 - e_1)(s_2 - k_1)} - \sqrt{(s_1 - e_1)(s_1 - k_1)(s_2 - e_2)(s_2 - e_3)(s_2 - k_2)}}{s_1 - s_2} \\
&= -C\frac{c_5}{c_{14}}\frac{\vartheta_{34}(v_1v_2)}{\vartheta_5(v_1v_2)},
\end{align*}
\begin{align*}
P_{15} &= \frac{\sqrt{(s_1 - e_2)(s_1 - e_3)(s_1 - k_1)(s_2 - e_1)(s_2 - k_2)} - \sqrt{(s_1 - e_1)(s_1 - k_2)(s_2 - e_2)(s_2 - e_3)(s_2 - k_1)}}{s_1 - s_2} \\
&= -C\frac{\vartheta_{13}(v_1v_2)}{\vartheta_5(v_1v_2)},
\end{align*}
\begin{align*}
P_{23} &= \frac{\sqrt{(s_1 - e_1)(s_1 - k_1)(s_1 - k_2)(s_2 - e_2)(s_2 - e_3)} - \sqrt{(s_1 - e_2)(s_1 - e_3)(s_2 - e_1)(s_2 - k_1)(s_2 - k_2)}}{s_1 - s_2} \\
&= C\frac{c_5}{c_4}\frac{\vartheta_{02}(v_1v_2)}{\vartheta_5(v_1v_2)},
\end{align*}
\begin{align*}
P_{24} &= \frac{\sqrt{(s_1 - e_1)(s_1 - e_3)(s_1 - k_2)(s_2 - e_2)(s_2 - k_1)} - \sqrt{(s_1 - e_2)(s_1 - k_1)(s_2 - e_1)(s_2 - e_3)(s_2 - k_2)}}{s_1 - s_2} \\
&= -C\frac{c_5}{c_0}\frac{\vartheta_{24}(v_1v_2)}{\vartheta_5(v_1v_2)},
\end{align*}
\begin{align*}
P_{25} &= \frac{\sqrt{(s_1 - e_1)(s_1 - e_3)(s_1 - k_1)(s_2 - e_2)(s_2 - k_2)} - \sqrt{(s_1 - e_2)(s_1 - k_2)(s_2 - e_1)(s_2 - e_3)(s_2 - k_1)}}{s_1 - s_2} \\
&= C\frac{c_5}{c_{23}}\frac{\vartheta_{12}(v_1v_2)}{\vartheta_5(v_1v_2)},
\end{align*}
\begin{align*}
P_{34} &= \frac{\sqrt{(s_1 - e_1)(s_1 - e_2)(s_1 - k_2)(s_2 - e_3)(s_2 - k_1)} - \sqrt{(s_1 - e_3)(s_1 - k_1)(s_2 - e_1)(s_2 - e_2)(s_2 - k_2)}}{s_1 - s_2} \\
&= -iC\frac{c_5}{c_2}\frac{\vartheta_{04}(v_1v_2)}{\vartheta_5(v_1v_2)},
\end{align*}
\begin{align*}
P_{35} &= \frac{\sqrt{(s_1 - e_1)(s_1 - e_2)(s_1 - k_1)(s_2 - e_3)(s_2 - k_2)} - \sqrt{(s_1 - e_3)(s_1 - k_2)(s_2 - e_1)(s_2 - e_1)(s_2 - k_2)}}{s_1 - s_2} \label{eq:P35} \\
&= -iC\frac{c_5}{c_{03}}\frac{\vartheta_{01}(v_1v_2)}{\vartheta_5(v_1v_2)}.
\end{align*}
\begin{align*}
P_{45} &= \frac{\sqrt{(s_1 - e_1)(s_1 - e_2)(s_1 - e_3)(s_2 - k_1)(s_2 - k_2)} - \sqrt{(s_1 - k_1)(s_1 - k_2)(s_2 - e_1)(s_2 - e_2)(s_2 - e_3)}}{s_1 - s_2} \\
&= -iC\frac{c_5}{c_{34}}\frac{\vartheta_{14}(v_1v_2)}{\vartheta_5(v_1v_2)}.
\end{align*}

Once the six quantities $p, q, r, \gamma, \gamma', \gamma''$ have been expressed in terms of the quotients $\frac{\vartheta_\alpha(v_1, v_2)}{\vartheta(v_1, v_2)}$, in which $v_1$, $v_2$ denote entire linear functions of time $t$, it remains to find the expressions for the six other cosines $\alpha, \alpha', \alpha'', \beta, \beta', \beta''$ as functions of time.

These quantities satisfy the differential equations
\begin{align*}
\frac{d\alpha}{dt} &= \alpha'r - \alpha''q, & \frac{d\beta}{dt} &= \beta'r - \beta''q, \\
\frac{d\alpha'}{dt} &= \alpha''p - \alpha r, & \frac{d\beta'}{dt} &= \beta''p - \beta r, \\
\frac{d\alpha''}{dt} &= \alpha q - \alpha'p, & \frac{d\beta''}{dt} &= \beta q - \beta'p.
\end{align*}

From these equations it follows that
\[
\frac{d(\alpha + \beta i)}{dt} = (\alpha' + \beta'i)r - (\alpha'' + \beta''i)q.
\]

Dividing by $\alpha + \beta i$ and remarking that
\begin{align*}
(\alpha' + \beta'i)(\alpha - \beta i) &= \alpha\alpha' + \beta\beta' + i(\alpha\beta' - \alpha'\beta) = -\gamma\gamma' + i\gamma'', \\
(\alpha'' + \beta''i)(\alpha - \beta i) &= \alpha\alpha'' + \beta\beta'' + i(\alpha\beta'' - \alpha''\beta) = -\gamma\gamma'' - i\gamma', \\
(\alpha + \beta i)(\alpha - \beta i) &= \alpha^2 + \beta^2 = 1 - \gamma^2,
\end{align*}
one finds
\begin{align*}
\frac{d}{dt}\log(\alpha + \beta i) &= r\frac{-\gamma\gamma' + i\gamma''}{1 - \gamma^2} + q\frac{\gamma\gamma'' + i\gamma'}{1 - \gamma^2} = \frac{-i\frac{d\gamma}{dt} + i(\gamma'q + \gamma''r)}{1 - \gamma^2} \\
&= \frac{1}{2}\left[\frac{\frac{d\gamma}{dt} + i(\gamma'q + \gamma''r)}{1 + \gamma} - \frac{\frac{d\gamma}{dt} - i(\gamma'q + \gamma''r)}{1 - \gamma}\right].
\end{align*}
In the same manner one finds
\begin{align*}
\frac{d}{dt}\log(\alpha' + \beta'i) &= \frac{1}{2}\left[\frac{\frac{d\gamma'}{dt} + i(\gamma p + \gamma''r)}{1 + \gamma'} - \frac{\frac{d\gamma'}{dt} - i(\gamma p + \gamma''r)}{1 - \gamma'}\right], \\
\frac{d}{dt}\log(\alpha'' + \beta''i) &= \frac{1}{2}\left[\frac{\frac{d\gamma''}{dt} + i(\gamma p + \gamma'q)}{1 + \gamma''} - \frac{\frac{d\gamma''}{dt} - i(\gamma p + \gamma'q)}{1 - \gamma''}\right].
\end{align*}

Substituting into the second of these equations the values of $p$, $q$, $r$, $\gamma$, $\gamma'$, $\gamma''$, expressed as functions of time, one sees that to obtain the values of $\alpha$, $\beta$, $\alpha'$, $\beta'$, $\alpha''$, $\beta''$ one must integrate rational functions of quotients $\frac{\vartheta_\alpha(t + c, c_1)}{\vartheta(t + c, c_1)}$.

One can demonstrate that each of the six quantities $\alpha$, $\alpha'$, $\alpha''$, $\beta$, $\beta'$, $\beta''$ is a uniform function of time, having only poles for finite values of time.
According to a well-known theorem, a uniform function $f(t)$ can be represented as a quotient of two convergent series for all finite values of $t$, provided the following condition is satisfied:
\[
\frac{d}{dt}\log f(t)
\]
can be expanded in the neighborhood of each finite value $t_0$ in a convergent series of the form
\[
\frac{d}{dt}\log f(t) = m(t - t_0)^{-1} + \mathfrak{P}(t - t_0)
\]
where $\mathfrak{P}(t - t_0)$ denotes an infinite series containing only terms with positive exponent, and $m$ is equal to zero or to a positive or negative integer.

Now, this is precisely the case for the right-hand sides of the equations which define
\[
\frac{d}{dt}\log(\alpha + \beta i), \quad \frac{d}{dt}\log(\alpha' + \beta'i), \quad \frac{d}{dt}\log(\alpha'' + \beta''i)
\]
as functions of time. If one writes there in place of $\gamma$, $\gamma'$, $\gamma''$, $p$, $q$, $r$ their values as functions of time and if, 
supposing for $t$ in the neighborhood of $t_0$, one expands them in a powers series of $t - t_0$, one sees immediately that terms with negative exponent can only enter in these expansions in the two following cases:

1) If the expansions of $p$, $q$, $r$, $\gamma$, $\gamma'$, $\gamma''$ contain negative powers of $(t - t_0)$.

2) If for $t = t_0$ one of the three quantities $\gamma$, $\gamma'$, $\gamma''$ is equal to $\pm 1$.

Having recourse to the differential equations (1) and writing $\tau$ in place of $t - t_0$, one sees that in the 1st case the expansions of $p$, $q$, $r$, $\gamma$, $\gamma'$, $\gamma''$ must have the following form:
\begin{align*}
p &= p_0\tau^{-1} + p_1 + p_2\tau + \ldots, & \gamma &= f_0\tau^{-2} + f_1\tau^{-1} + f_2 + \ldots, \\
q &= q_0\tau^{-1} + q_1 + q_2\tau + \ldots, & \gamma' &= g_0\tau^{-2} + g_1\tau^{-1} + g_2 + \ldots, \\
r &= r_0\tau^{-1} + r_1 + r_2\tau + \ldots, & \gamma'' &= h_0\tau^{-2} + h_1\tau^{-1} + h_2 + \ldots,
\end{align*}
where the coefficients $p_0$, $q_0$, $r_0$, $f_0$, $g_0$, $h_0$ can have these two different systems of values (see § 1):
\begin{align*}
\text{I.} \quad r_0 &= 0, & h_0 &= \pm i\frac{4}{c_0}, \\
q_0 &= \pm i2, & g_0 &= 0, \\
p_0 &= 0, & f_0 &= -\frac{4}{c_0}, \\
\text{II.} \quad q_0 &= \pm ip_0, & g_0 &= \pm if_0 = -\frac{r_0}{c_0} = \pm \frac{2i}{c_0}, \\
r_0 &= \pm 2i, & h_0 &= 0.
\end{align*}

Substituting these expansions into the right hand side of the equation
\begin{align*}
\frac{d}{dt}\log(\alpha'' + \beta''i) &= \frac{\gamma''\frac{d\gamma''}{dt} - i(\gamma p + \gamma'q)}{\gamma''^2 - 1} \\
&= \frac{1}{2}\left[\frac{\frac{d\gamma''}{dt} - i(\gamma p + \gamma'q)}{\gamma'' - 1} + \frac{\frac{d\gamma''}{dt} + i(\gamma p + \gamma'q)}{\gamma'' + 1}\right],
\end{align*}

one sees that one has in case I
\[
\frac{d}{dt}\log(\alpha'' + \beta''i) = -2\tau^{-1} + \mathfrak{P}(\tau)
\]
and in case II
\[
\frac{d}{dt}\log(\alpha'' + \beta''i) = -\tau^{-1} + \mathfrak{P}(\tau).
\]

Let us now examine more closely the expansion of $\frac{d}{dt}\log(\alpha'' + \beta''i)$ in the neighborhood of a value $t = t_0$ for which $\gamma'' = \pm 1$.

In general, if one sets
\begin{align*}
p &= p_0 + p_1t + \ldots, & \gamma &= f_0 + f_1t + \ldots, \\
q &= q_0 + q_1t + \ldots, & \gamma' &= g_0 + g_1t + \ldots, \\
r &= r_0 + r_1t + \ldots, & \gamma'' &= h_0 + h_1t + \ldots,
\end{align*}
the coefficients $p_1$, $q_1$, $r_1$, $f_1$, $g_1$, $h_1$ are defined (by virtue of the differential equations (1), § 2) as functions of $p_0$, $q_0$, $r_0$, $f_0$, $g_0$, $h_0$ by the following equations
\begin{align*}
2p_1 &= q_0r_0, & f_1 &= r_0g_0 - q_0h_0, \\
2q_1 &= -p_0r_0 - c_0h_0, & g_1 &= p_0h_0 - r_0f_0, \\
r_1 &= c_0g_0, & h_1 &= q_0f_0 - p_0g_0.
\end{align*}

The three quantities $\gamma$, $\gamma'$, $\gamma''$ are linked by the equation
\[
\gamma^2 + \gamma'^2 + \gamma''^2 = 1.
\]

Their initial values $f_0$, $g_0$, $h_0$ are therefore also subject to the condition
\[
f_0^2 + g_0^2 + h_0^2 = 1.
\]

Let us denote by $\varepsilon_1$ and by $\varepsilon_2$ two quantities which satisfy the equation
\[
\varepsilon_1^2 = 1, \quad \varepsilon_2^2 = 1.
\]

If we set $h_0 = \varepsilon_1$ it is therefore necessary that one has at the same time $g_0 = \varepsilon_2if_0$ and one then finds
\begin{align*}
h_1 &= -\varepsilon_2if_0(p_0 + \varepsilon_2iq_0), \\
p_0f_0 + q_0g_0 &= f_0(p_0 + \varepsilon_2iq_0) = \varepsilon_2ih_1.
\end{align*}

The expansion of the term
\[
\frac{1}{2}\frac{\frac{d\gamma''}{dt} - i\varepsilon_1(p\gamma + q\gamma')}{\gamma'' - \varepsilon_1}
\]
which alone can contain negative powers of $\tau$, therefore has the form
\[
\frac{1 + \varepsilon_1\varepsilon_2}{2}\tau^{-1} + \mathfrak{P}(\tau).
\]

If $\varepsilon_1$ and $\varepsilon_2$ have the same sign the coefficient of $\tau^{-1}$ is equal to $+1$.

If $\varepsilon_1$ and $\varepsilon_2$ have contrary signs this coefficient is null.

One therefore sees that in all cases, the expansion of
\[
\frac{d}{dt}\log(\alpha'' + \beta''i)
\]
has the required form
\[
\frac{d}{dt}\log(\alpha'' + \beta''i) = m(t - t_0)^{-1} + \mathfrak{P}(t - t_0),
\]
from which it results that $\alpha'' + \beta''i$ is a uniform function of time, which can be put in the form of a quotient of two always convergent series.

One arrives at the same result concerning $\alpha + \beta i$ and $\alpha' + \beta'i$.

I have also ascertained that these quantities can be expressed as rational functions of quantities of the form
\[
\frac{\vartheta_\alpha(u_1 + v_1, u_2 + v_2)}{\vartheta(u_1, u_2)}e^{u_3},
\]
$\alpha$ being the index of one of 16 functions $\vartheta(u_1, u_2)$, $u_1$, $u_2$, $u_3$ denoting entire, linear functions of time, and $v_1$, $v_2$ denoting imaginary constants.

But because of the considerable complexity of the calculations, I have not yet succeeded in deriving these formulas in their final form.

\section*{§ 8.}
\setcounter{equation}{0}

It is interesting to realise a model for the case of rotation of a rigid body about a fixed point where the conditions of the case just studied are fulfilled.

In the differential equations (1) the constants $A$, $B$, $C$ denote the three principal moments of inertia with respect to the fixed point.
Let now $A_1$, $B_1$, $C_1$ be the three principal moments of inertia of the same body with respect to its center of gravity and let us denote by $x$, $y$, $z$ the coordinates of a point in a system of coordinate axes whose origin is at the center of gravity and whose directions coincide with the three principal axes of inertia passing through this point.
One then has, denoting by $\mu$ the density of the considered body at the point whose coordinates are $x$, $y$, $z$,
\begin{align*}
A_1 &= \iiint\mu(y^2 + z^2)\,dxdydz, \qquad B_1 = \iiint\mu(z^2 + x^2)\,dxdydz, \qquad C_1 = \iiint\mu(x^2 + y^2)\,dxdydz, 
\end{align*}
\begin{align*}
0 &= \iiint\mu x\,dxdydz, \qquad 0 = \iiint\mu y\,dxdydz, \qquad 0 = \iiint\mu z\,dxdydz, \\
0 &= \iiint\mu yz\,dxdydz, \qquad 0 = \iiint\mu zx\,dxdydz, \qquad 0 = \iiint\mu xy\,dxdydz.
\end{align*}
These triple integrals must be extended over the entire interior of the considered body.

Let $a$, $b$, $c$ be the coordinates of a point $O$. Denoting by $\xi$, $\eta$, $\zeta$ the coordinates of a point of space with respect to axes parallel to the preceding ones, but whose origin coincides with point $O$, one has
\[
\xi = x - a, \quad \eta = y - b, \quad \zeta = z - c.
\]

The equation of the ellipsoid of inertia whose center is at $O$ is, as is well known,
\[
1 = A'_1\xi^2 + B'_1\eta^2 + C'_1\zeta^2 - 2D'_1\eta\zeta - 2E'_1\zeta\xi - 2F'_1\xi\eta = \varphi(\xi, \eta, \zeta),
\]
where the coefficients $A'_1$, $B'_1$, \ldots are defined by the equations
\begin{align*}
A'_1 &= \iiint\mu(\eta^2 + \zeta^2)\,d\xi d\eta d\zeta = A_1 + M(b^2 + c^2), \\
B'_1 &= \iiint\mu(\zeta^2 + \xi^2)\,d\xi d\eta d\zeta = B_1 + M(c^2 + a^2), \\
C'_1 &= \iiint\mu(\xi^2 + \eta^2)\,d\xi d\eta d\zeta = C_1 + M(a^2 + b^2), \\
D'_1 &= \iiint\mu\eta\zeta\, d\xi d\eta d\zeta = Mbc, \\
E'_1 &= \iiint\mu\zeta\xi\, d\xi d\eta d\zeta = Mca, \\
F'_1 &= \iiint\mu\xi\eta\, d\xi d\eta d\zeta = Mab,
\end{align*}
and where $M$ denotes the total mass of the body,
\[
M = \iiint\mu\, d\xi d\eta d\zeta.
\]

Let now $u$, $v$, $w$ be the coordinates of a point in space described relative to axes with origin at point $O$, but with directions coinciding with those of the principal axes of the ellipsoid of inertia, corresponding to this point. One then has
\begin{align*}
u &= \alpha\xi + \beta\eta + \gamma\zeta, \\
v &= \alpha_1\xi + \beta_1\eta + \gamma_1\zeta, \\
w &= \alpha_2\xi + \beta_2\eta + \gamma_2\zeta,
\end{align*}
\begin{align*}
\xi^2 + \eta^2 + \zeta^2 &= u^2 + v^2 + w^2, \\
\varphi(\xi, \eta, \zeta) &= \lambda u^2 + \mu v^2 + \nu w^2,
\end{align*}
where $\lambda$, $\mu$, $\nu$ are positive constants.
Moreover, denoting by $u_0$, $v_0$, $w_0$ the coordinates of the center of gravity in this new system, one has
\begin{align*}
-u_0 &= \alpha a + \beta b + \gamma c, \\
-v_0 &= \alpha_1 a + \beta_1 b + \gamma_1 c, \\
-w_0 &= \alpha_2 a + \beta_2 b + \gamma_2 c.
\end{align*}

If we now impart to the considered body a motion of rotation about the fixed point $O$, the conditions of the case studied by us will be fulfilled, if one has
\[
\lambda = \mu = 2\nu,
\]
\[
-w_0 = \alpha_2 a + \beta_2 b + \gamma_2 c = 0.
\]
One must therefore be able to satisfy the following equations
\begin{align*}
\varphi(\xi, \eta, \zeta) &= \nu\left[\left(2u^2+v^2\right)+w^2\right], & a \alpha_2  +  b \beta_2 +  c \gamma_2 &= 0, \\
\xi^2 + \eta^2 + \zeta^2 &= u^2 + v^2 + w^2, & \alpha_2^2 + \beta_2^2 + \gamma_2^2 &= 1.
\end{align*}

Let us see if one can choose the point $O$ in such a manner that these equations are fulfilled.

From the first two of these equations, it follows that
\[
\varphi(\xi, \eta, \zeta) - 2\nu(\xi^2 + \eta^2 + \zeta^2) = -\nu w^2.
\]
This last equation must be satisfied, if when one writes for $w$ its value
\[
w = \alpha_2\xi + \beta_2\eta + \gamma_2\zeta,
\]
one obtains, by equating to zero the coefficients of each term in this identity, the following equations:
\begin{align*}
A'_1 - 2\nu &= A_1 + M(b^2 + c^2) - 2\nu = -\nu\alpha_2^2, \\
B'_1 - 2\nu &= B_1 + M(c^2 + a^2) - 2\nu = -\nu\beta_2^2, \\
C'_1 - 2\nu &= C_1 + M(a^2 + b^2) - 2\nu = -\nu\gamma_2^2, \\
D'_1 &= Mbc = \nu\beta_2\gamma_2, \\
E'_1 &= Mca = \nu\gamma_2\alpha_2, \\
F'_1 &= Mab = \nu\alpha_2\beta_2.
\end{align*}
If none of the three constants $a$, $b$, $c$ (the coordinates of point $O$) were zero it would follow from the last three equations that
\[
M^3a^2b^2c^2 = \nu^3\alpha_2^2\beta_2^2\gamma_2^2,
\]
\begin{align*}
M^{\frac{1}{2}}a &= \nu^{\frac{1}{2}}\alpha_2, \\
M^{\frac{1}{2}}b &= \nu^{\frac{1}{2}}\beta_2, \\
M^{\frac{1}{2}}c &= \nu^{\frac{1}{2}}\gamma_2,
\end{align*}
which is evidently impossible as a result of the equations
\[
\alpha_2^2 + \beta_2^2 + \gamma_2^2 = 1, \quad a \alpha_2  + b \beta_2  + c \gamma_2  = 0,
\]
with neither $M$ nor $\nu$ equal to 0.

If we suppose $c = 0$, but $a$ and $b$ different from 0, we must set $\gamma_2 = 0$, and consequently
\[
Mab = \nu\alpha_2\beta_2, \quad \alpha_2^2 + \beta_2^2 = 1, \quad a \alpha_2 + b \beta_2 = 0,
\]
from which it follows that
\[
\alpha_2 = \frac{b}{\sqrt{a^2 + b^2}}, \quad \beta_2 = \frac{-a}{\sqrt{a^2 + b^2}},
\]
(the sign of $\sqrt{a^2 + b^2}$ being determined arbitrarily) leading to
\[
M = \frac{-\nu}{a^2 + b^2},
\]
which is impossible to satisfy since $M$ and $\nu$ are both positive quantities.

It is therefore necessary to suppose $b = 0$, $c = 0$. One then has
\begin{align*}
A'_1 &= A_1, & B'_1 &= B_1 + Ma^2, & C'_1 &= C_1 + Ma^2, \\
D'_1 &= 0, & E'_1 &= 0, & F'_1 &= 0,
\end{align*}
\[
\varphi(\xi, \eta, \zeta) = A_1\xi^2 + (B_1 + Ma^2)\eta^2 + (C_1 + Ma^2)\zeta^2.
\]

If the three principal axes of inertia $A_1$, $B_1$, $C_1$ relative to the center of gravity of the considered body satisfy the equation
\[
A_1 = 2(B_1 - C_1),
\]
one will be able to satisfy all the conditions supposed by us, by taking
\[
a^2 = \frac{A_1 - B_1}{M},
\]
for one has in this case
\[
B_1 + Ma^2 = A_1, \qquad C_1 + Ma^2 = C_1 + A_1 - B_1 = \frac{1}{2}A_1,
\]
and consequently
\[
\varphi(\xi, \eta, \zeta) = A_1\left(\xi^2 + \eta^2 + \frac{1}{2}\zeta^2\right).
\]

Let us remark only that for $a$ to be real, it is necessary and sufficient that one has $B_1 > 2C_1$.

One therefore sees, according to this calculation, that it is possible to realise mechanically all the conditions of the problem that I have just studied.

\vspace{1in}
\begin{center}
\raisebox{0.5ex}{\rule{1in}{0.5pt}}\quad$\bullet$\quad$\diamond$\quad$\bullet$\quad\raisebox{0.5ex}{\rule{1in}{0.5pt}}
\end{center}

\newpage
\section*{Translator's Notes}

\noindent\hypertarget{TN1}{}\textbf{\hyperlink{TN1-ref}{[TN1]} Equation numbering.} The author numbers equation \eqref{eq:14} as (14), though it could logically have been numbered (10) given that the previous numbered equation was (9). The original numbering has been retained in this translation.

\noindent\hypertarget{TN2}{}\textbf{\hyperlink{TN2-ref}{[TN2]} Possible typographical error in equation 14.} In equation \eqref{eq:14}, the final term of the expression for $h_2'''$ contains $(e_1^2 - e_3^2)$ in the original text. Based on the symmetry pattern observed in the preceding expressions for $h_2'$ and $h_2''$, this may be a typographical error for $(e_1^2 - e_2^2)$. However, the original notation has been preserved in this translation.

\noindent\hypertarget{TN3}{}\textbf{\hyperlink{TN3-ref}{[TN3]} Possible typographical error in expression for $P_{35}$.} The second square root in the numerator of the expression for $P_{35}$ contains the factor $(s_2 - e_1)$ repeated twice. By analogy with the patterns in the expressions for $P_{23}$, $P_{24}$, $P_{25}$, and $P_{34}$, one of these factors should likely be $(s_2 - e_2)$. However, the original notation has been preserved in this translation.

\vspace{1em}
\noindent\textbf{Original paper.} Kowalevski, S. (1889). Sur le problème de la rotation d'un corps solide autour d'un point fixe. \textit{Acta Mathematica}, 12, 177--232. Available at: \url{https://doi.org/10.1007/BF02592182}

\vspace{1em}

\noindent{Translated by G. Hesketh 06/03/2026 and dedicated to Mich\`ele Audin: \textit{``And if we do not try\ldots\ I try.''}}
\end{document}